\documentclass[twoside,12pt]{article}
\usepackage{amsmath}
\usepackage{amssymb}
\usepackage{cases}
\usepackage{cite}
\usepackage[usenames,dvipsnames]{color}
\definecolor{rosso}{rgb}{0.8,0,0}
\def\takeshi{\color{red}}
\def\pier{\color{red}}
\def\revis{\color{red}}
\def\colli{\color{red}}
\def\fukao{\color{red}}

\let\takeshi\relax
\let\pier\relax
\let\revis\relax
\let\colli\relax
\let\fukao\relax

\renewcommand{\theequation}{\thesection.\arabic{equation}}

\title{{C}ahn--{H}illiard equation on the boundary with bulk condition of {A}llen--{C}ahn type}

\author{Pierluigi Colli\\
Dipartimento di Matematica, Universit\`a degli Studi di Pavia\\
Via Ferrata~1, 27100 Pavia, Italy\\
E-mail: \texttt{pierluigi.colli@unipv.it}\\
\and \\ Takeshi Fukao\\
Department of Mathematics, Faculty of Education\\
Kyoto University of Education\\
1~Fujinomori, Fukakusa, Fushimi-ku, Kyoto~612-8522 Japan\\
E-mail: \texttt{fukao@kyokyo-u.ac.jp}}
 
\date{}

\newcommand\testopari{\sc Pierluigi Colli and Takeshi Fukao}
\newcommand\testodispari{\sc {C}ahn--{H}illiard equation on the boundary with bulk condition}
\begin{document}

\maketitle

\begin{abstract}
The well-posedness 
for a system of partial differential equations and dynamic boundary conditions is discussed. 
This system is 
a sort of transmission problem between the 
dynamics in the bulk $\Omega $ and on the boundary $\Gamma$. 
The {P}oisson equation for the chemical potential, 
the {A}llen--{C}ahn equation for the order parameter in the bulk $\Omega$ 
are considered 
as auxiliary conditions for
solving the {C}ahn--{H}illiard equation on the boundary $\Gamma$. 
Recently the well-posedness for the equation 
and dynamic boundary condition, both of {C}ahn--{H}illiard type, was 
discussed. 
Based on this result,
 the existence of the solution and its continuous dependence on the data are 
proved.

\vspace{2mm}

\vspace{2mm}

\end{abstract}

\section{Introduction}
\setcounter{equation}{0}

In this paper, 
we treat the {C}ahn--{H}illiard equation \cite{CH58} on the boundary of some 
bounded smooth domain. 
Let $0<T<+\infty$ be some fixed time and let $\Omega \subset \mathbb{R}^{d}$, 
$d=2$ or $3$, be a bounded domain occupied by a material and its boundary 
$\Gamma$ of $\Omega $ is supposed to be smooth enough.  
We start from the following equations of 
{C}ahn--{H}illiard type on the boundary $\Gamma$: 
\begin{gather} 
	\partial_t u_\Gamma -\Delta_\Gamma \mu_\Gamma  = -\partial _{\boldsymbol{\nu }} \mu 
	\quad \mbox{on }\Sigma:=\Gamma \times (0,T), 
	\label{CHB1}
	\\
	\mu_\Gamma = -\Delta_\Gamma u_\Gamma  + {\mathcal W}_\Gamma'(u_\Gamma)-f_\Gamma +\partial _{\boldsymbol{\nu }} u
	\quad \mbox{on }\Sigma,
	\label{CHB2}
\end{gather}
where $\partial_t$ denotes the partial 
derivative with respect to time, and 
$\Delta _{\Gamma }$ denotes the {L}aplace--{B}eltrami operator
on $\Gamma $ (see, e.g., \cite[Chapter~3]{Gri09}). 
Here, the
unknowns $u_\Gamma$ and $\mu_\Gamma:\Sigma \to \mathbb{R}$ 
stand for the order parameter and
the chemical potential, respectively. 
In the right hand sides of \eqref{CHB1} and \eqref{CHB2}, 
the outward normal derivative $\partial_{\boldsymbol{\nu }}$ on $\Gamma $ 
acts on functions 
$\mu, u:Q:=(0,T)\times \Omega \to \mathbb{R}$ that 
satisfy the following trace conditions
\begin{equation}
	\mu _{|_\Sigma } = \mu _\Gamma, 
	\quad 
	u_{|_\Sigma } =u_\Gamma 
	\quad \mbox{on } \Sigma,  
	\label{traceor}
\end{equation}
where $\mu  _{|_\Sigma}$ and $u _{|_\Sigma}$ are the traces of 
$\mu $ and $u$ on $\Sigma$. 
Moreover, these functions $\mu $ and $u$ solve the following 
equations in the bulk $\Omega$:
\begin{gather} 
	- \Delta \mu =0
	\quad \mbox{in } Q,
	\label{P}
	\\
	\tau \partial _t u - \Delta u +{\mathcal W}' (u) = f
	\quad \mbox{in } Q,
	\label{AC}
\end{gather}
where $\tau >0$ is a positive constant and 
$\Delta$ denotes the {L}aplacian.

Note that the nonlinear terms 
${\mathcal W}'_\Gamma$ and ${\mathcal W}'$ are the derivatives 
of the functions ${\mathcal W}_\Gamma$ and ${\mathcal W}$ usually referred as double-well potentials,
with two minima and a local unstable maximum in between. 
The prototype model 
is provided by ${\mathcal W}_\Gamma (r)={\mathcal W}(r)=(1/4)(r^2-1)^2$ 
so that ${\mathcal W}'_\Gamma (r)={\mathcal W} '(r)=r^3-r$, $r\in \mathbb{R}$, is the sum of 
an increasing function with a power growth and another smooth 
function which breaks the monotonicity properties of 
the former and is related to the non-convex part 
of the potential ${\mathcal W}_\Gamma $ or ${\mathcal W}$.

Therefore, we can say that the system \eqref{CHB1}--\eqref{CHB2} yields 
the {C}ahn--{H}illiard equation on a smooth manifold $\Gamma$, 
and equations \eqref{P} 
and \eqref{AC} reduce to the {P}oisson equation for $\mu $ and 
the {A}llen--{C}ahn equation for $u$ in the bulk $\Omega$, as auxiliary conditions for
solving \eqref{CHB1}--\eqref{CHB2}. In other word, \eqref{CHB1}--\eqref{AC} is 
a sort of transmission problem between the 
dynamics in the bulk $\Omega $ and the one on the boundary $\Gamma$. 
With the following initial conditions 
\begin{equation} 
	u_\Gamma  (0)=u_{0\Gamma}
	\quad \mbox{on }\Gamma,
	\quad 
	u(0)=u_0
	\quad 
	\mbox{in }\Omega, 
	\label{IC}
\end{equation} 
the problem 
\eqref{CHB1}--\eqref{IC} becomes the initial value problem of 
a {P}oisson--{A}llen--{C}ahn system with the dynamic 
boundary condition of {C}ahn--{H}illiard type, named by (P). Indeed, the 
interaction between $\mu $ and $u$ appears only on \eqref{CHB1}--\eqref{CHB2}, 
whereas \eqref{P} and \eqref{AC} are independent equations. 
As a remark, if $\tau =0$ then the problem (P) turns out to be a \emph{quasi-static system}. 
From \eqref{CHB1}, \eqref{P} and \eqref{IC}, we easily see that the mass conservation on the 
boundary holds as follows:
\begin{equation}
	\int_{\Gamma }^{} u_{\Gamma }(t) d\Gamma 
	= \int_{\Gamma }^{} u_{0\Gamma } d\Gamma
	\quad \mbox{for all }t \in [0,T].
	\label{masscon}
\end{equation}

Let us mention some related results: 
the {\revis papers \cite{Esc94, FQ97, FIK12, FIK15, GM14} consider} 
some quasi-static systems with 
dynamic boundary conditions (see also \cite{CF15, CF15a}); 
the contributions{\pier \cite{CGM13, CF15b, CGS14, CGS17, Gal06, GW08, GMS11, LW17}} set 
a {C}ahn--{H}illiard equation on the boundary as the dynamic boundary condition;
a more complicated system of {C}ahn--{H}illiard type on the boundary 
with a mass conservation condition is investigated 
in \cite{GKRR16}. 
Especially, in this paper we will exploit the previous result in \cite{CF15b}, on which 
equations \eqref{CHB1}--\eqref{AC} were replaced by 
\begin{gather}
	\varepsilon \partial _t u -\Delta \mu  =0 
	\quad \mbox{in } Q, 
	\label{Pe}\\
	\varepsilon \mu 
	=\tau \partial_t u - \Delta u 
	+ {\mathcal W}'(u) -f 
	\quad \mbox{in } Q,
	\label{ACe}
	\\
	u_{|_\Sigma } =u_\Gamma, \quad 
	\mu _{|_\Sigma } = \mu _\Gamma 
	\quad \mbox{on } \Sigma, 
	\label{traceore}
	\\
	\partial_t u_\Gamma +\partial _{\boldsymbol{\nu }} \mu-\Delta_\Gamma \mu_\Gamma  = 0
	\quad \mbox{on }\Sigma, 
	\label{CHB1e}
	\\
	\mu_\Gamma = \partial _{\boldsymbol{\nu }} u-\Delta_\Gamma u_\Gamma  + {\mathcal W}_\Gamma'(u_\Gamma)-f_\Gamma 
	\quad \mbox{on }\Sigma,
	\label{CHB2e}
\end{gather}
where $\varepsilon >0$. Then from the system 
\eqref{Pe}--\eqref{CHB2e} with \eqref{IC}, we 
obtain 
the 
following total mass conservation:
\begin{equation}
	\varepsilon \int_{\Omega }^{} u(t) dx + \int_{\Gamma }^{} u_{\Gamma }(t) d\Gamma 
	= \varepsilon \int_{\Omega }^{} u_0  dx+\int_{\Gamma }^{} u_{0\Gamma } d\Gamma
	\quad \mbox{for all }t \in [0,T].
	\label{totalmass}
\end{equation}
Our essential idea of the existence proof is to be able to pass to the limit as $\varepsilon \searrow  0$ in 
the system 
\eqref{Pe}--\eqref{CHB2e} with \eqref{IC}. 
To be more precise about our arguments, let us introduce 
a brief outline of the paper along with 
a short description of the various items.

In Section~2, we present the main results of 
the well-posedness of the system \eqref{CHB1}--\eqref{IC}. 
A solution to the problem (P) is suitably defined. 
The main theorems are concerned with the existence of the solution (Theorem~2.1) 
and the continuous dependence on the given data (Theorem~2.2), 
the second theorem entailing the uniqueness property.

In Section~3, we consider the approximate problem for (P) with 
two approximate parameters $\varepsilon $ and $\lambda $ 
by substituting the maximal monotone graphs with their 
{Y}osida regularizations
in terms of the parameter $\lambda $. 
Moreover, we obtain the uniform estimates with suitable growth order. 
Here, we can apply the results obtained in \cite{CF15b}.

In Section~4, we prove the existence result. 
The proof is split in several steps. 
In the first step, we obtain the uniform estimates with respect to 
$\varepsilon $. Then combining with the previous estimates of Section~3,
 we can pass to the limit as $\varepsilon \searrow  0$. 
In the second step, we improve suitable estimates in order to make them 
independent of $\lambda$. 
Then, we can pass to the limit as $\lambda \searrow 0$ 
and conclude the existence proof. 
The last part of this section is devoted to the proof of the continuous dependence.

Finally, the last section contains an Appendix 
in which the approximate problem for (P) and 
some auxiliary results are discussed.



\section{Main results}
\setcounter{equation}{0}

In this section, our main result is stated. 
At first, we give our target system {\rm (P)} of equations and conditions 
as follows: for any fixed constant $\tau >0$, we have 
\begin{gather} 
	-\Delta \mu = 0 
	\quad \mbox{a.e.\ in }Q,
	\label{QS1}
\\
	\tau \partial _t u-\Delta u + \xi + \pi (u) =f, 
	\quad \xi \in \beta (u)
	\quad \mbox{a.e.\ in }Q,
	\label{QS2}
\\
	u_\Gamma =u _{|_\Sigma }, \quad \mu _\Gamma =\mu _{|_\Sigma},
	\quad 
	{\partial_t u_\Gamma }+\partial_{\boldsymbol{\nu }} \mu -\Delta _\Gamma \mu _\Gamma =0
	\quad \mbox{a.e.\ on }\Sigma,
	\label{QS3}
\\ 
	\mu _\Gamma =\partial_{\boldsymbol{\nu }} u - \Delta _\Gamma u_\Gamma +\xi _\Gamma + \pi _\Gamma (u_\Gamma )-f_\Gamma, 
	\quad \xi_\Gamma  \in \beta _\Gamma (u_\Gamma )
	\quad \mbox{a.e.\ on }\Sigma,
	\label{QS4}
\\
	u(0)=u_0 
	\quad 
	\mbox{a.e.\ in }\Omega, \quad 
	u_\Gamma  (0)=u_{0\Gamma}
	\quad \mbox{a.e.\ on }\Gamma,
	\label{QS5}
\end{gather} 
where $f:Q \to \mathbb{R}$, $f_\Gamma :\Sigma \to \mathbb{R}$ are given sources, 
$u_0:\Omega \to \mathbb{R}$, $u_{0\Gamma} :\Gamma \to \mathbb{R}$ are known initial data; 
$\beta$ stands for the subdifferential of the convex part $\widehat{\beta }$ and 
$\pi $ stands for the derivative of the concave perturbation $\widehat{\pi}$ of a 
double well potential ${\mathcal W}(r)=\widehat{\beta }(r)+\widehat{\pi}(r)$ for all 
$r \in \mathbb{R}$. 
The same setting holds for $\beta _\Gamma $ and $\pi _\Gamma $.

Typical examples of the nonlinearities {\colli $\beta$, $\pi $ are given by:
\begin{itemize}
	\item $\beta(r)=r^3$, $\pi (r)=-r$ 
	for all $r \in \mathbb{R}$ with $D(\beta )=\mathbb{R}$ 
	for the prototype double well potential
	\begin{equation*} 
	{\mathcal W}(r)=\frac{1}{4}(r^2-1)^2;
	\end{equation*} 
	\item $\beta(r)=\ln((1+r)/(1-r))$, $\pi (r)=-2cr$ 
	for all $r \in D(\beta )$ with $D(\beta )=(-1,1)$ 
	for the logarithmic double well potential 
	\begin{equation*} 
	{\mathcal W}(r)= \bigl( (1+r) \ln(1+r)+(1-r)\ln(1-r) \bigr)-cr^2,
	\end{equation*} 
	where $c>0$ is a large constant which breaks convexity;
	\item $\beta(r)=\partial I_{[-1,1]}(r)$, $\pi (r)=-r$ 
	for all $r \in D(\beta )$ with $D(\beta )=[-1,1]$ 
	for the singular potential 
	\begin{equation*}
	{\mathcal W}(r)=I_{[-1,1]}(r)- \frac{1}{2}r^2,
	\end{equation*} 
	where $I_{[-1,1]}$ is the indicator function on $[-1,1]$. 
\end{itemize}
Similar choices can be considered for $\beta_\Gamma, \pi _\Gamma$ and the related potential 
${\mathcal W}_\Gamma $. What is important in our approach is that the potential on the boundary should dominate the potential in the bulk, that is, we prescribe a compatibility condition 
between $\beta$ and $\beta_\Gamma $ (see the later assumption (A5)) that forces the growth of $\beta $ to be controlled by the growth of $\beta_\Gamma$. A similar approach was taken in previous analyses, see \cite{CC13, CF15, CF15a, CF15b,CGS17, LW17}.}

As a remark, $\tau >0$ plays the role of a viscous parameter: indeed,
if $\tau =0$, then equations \eqref{QS1} and \eqref{QS2} become the 
stationary problem in $Q$, namely the quasi-static system. 
{\takeshi {\colli A natural question arises whether one can investigate also the case 
$\tau =0$ or, in our framework, also study the} singular limit {\colli as} $\tau \searrow  0${\colli : in our opinion, 
this is not a trivial question and {\pier deserves some attention and efforts. For}
the moment, we can just highlight it as open problem.}}

\subsection{Definition of the solution}

We treat the problem {\rm (P)} by a system of weak formulations. 
To do so, we introduce the
spaces $H:=L^2(\Omega )$, 
$H_\Gamma :=L^2(\Gamma )$, $V:=H^1(\Omega )$, $V_\Gamma :=H^1(\Gamma )$, 
$W:=H^2(\Omega)$, $W_\Gamma :=H^2(\Gamma)$
with usual norms and inner products, 
denote them by  
$| \cdot |_{H}$, $| \cdot |_{H_\Gamma}$ and $(\cdot,\cdot )_{H}$, $(\cdot,\cdot )_{H_\Gamma}$, and so on. {\pier Concerning these inner products and norms, we use the same notation for scalar and vectorial functions (e.g., typically gradients).}

Moreover, we put $\boldsymbol{H}:=H \times H_\Gamma $ and 
$\boldsymbol{V}:=\{ (z,z_\Gamma ) \in V \times V_\Gamma : z_\Gamma =z_{|_\Gamma } \ \mbox{a.e.\ on } \Gamma \}$.
Then, $\boldsymbol{H}$ and $\boldsymbol{V}$ are {H}ilbert spaces 
with the inner products 
\begin{align*}
	(\boldsymbol{u},\boldsymbol{z})_{\boldsymbol{H}}
	:=(u,z)_{H} + (u_\Gamma ,z_\Gamma )_{H_\Gamma } \quad 
	& \mbox{for all } \boldsymbol{u}:=(u,u_{\Gamma }), 
	\boldsymbol{z}:=(z,z_{\Gamma }) \in \boldsymbol{H},\\
	(\boldsymbol{u},\boldsymbol{z})_{\boldsymbol{V}}
	:=(u,z)_{V} + (u_\Gamma ,z_\Gamma )_{V_\Gamma } \quad 
	& \mbox{for all } \boldsymbol{u}
	:=(u,u_{\Gamma }), \boldsymbol{z}:=(z,z_{\Gamma }) 
	\in \boldsymbol{V}
\end{align*}
and related norms. As a remark, 
if 
$\boldsymbol{z}:=(z,z_{\Gamma}) \in \boldsymbol{V}$ then $z_{\Gamma }$ 
is the trace $z_{|_\Gamma }$ of $z$ on $\Gamma$, 
while if 
$\boldsymbol{z}:=(z,z_{\Gamma}) \in \boldsymbol{H}$ then 
$z \in H$ and $z_{\Gamma } \in H_{\Gamma }$ are independent. 
Hereafter, we use the notation of a bold letter 
like $\boldsymbol{z}$ to denote 
the pair which corresponds to the letter, that is $(z,z_\Gamma )$ for $\boldsymbol{z}$.

It is easy to see that the problem {\rm (P)} has a structure of volume conservation on the boundary 
$\Gamma$. Indeed, integrating the last equation in \eqref{QS3} on $\Sigma$, and using \eqref{QS1} and \eqref{QS5} we obtain
\begin{equation}
	\int_{\Gamma }^{} u_{\Gamma }(t) d\Gamma 
	= \int_{\Gamma}^{} u_{0\Gamma } d\Gamma
	\quad \mbox{for all }t \in [0,T];
	\label{mass}
\end{equation}
hereafter, we put 
\begin{equation}
	m_\Gamma:=\frac{1}{|\Gamma |}\int_{\Gamma }^{}u_{0\Gamma } d\Gamma,
	\label{mgamma}
\end{equation}
where $|\Gamma |:=\int_{\Gamma }^{}1 d\Gamma$. 
The space $\boldsymbol{V}^*$ denotes the dual of $\boldsymbol{V}$, 
and $\langle \cdot ,\cdot \rangle _{\boldsymbol{V}^*, \boldsymbol{V}}$ denotes
the duality pairing between $\boldsymbol{V}^*$ and 
$\boldsymbol{V}$. Moreover, it is understood that $\boldsymbol{H}$ is embedded in 
$\boldsymbol{V}^*$ in the usual way, i.e., 
$\langle \boldsymbol{u}, \boldsymbol{z} \rangle_{\boldsymbol{V}^*,\boldsymbol{V}}
= (\boldsymbol{u},\boldsymbol{z})_{\boldsymbol{H}}$ for all 
$\boldsymbol{u} \in \boldsymbol{H}$, $\boldsymbol{z} \in \boldsymbol{V}$. 
Then, we obtain 
$\boldsymbol{V}
\mathop{\hookrightarrow} \mathop{\hookrightarrow}
\boldsymbol{H}
\mathop{\hookrightarrow} \mathop{\hookrightarrow}
\boldsymbol{V}^*$, where 
``$\mathop{\hookrightarrow} \mathop{\hookrightarrow} $'' stands for 
the dense and compact embedding, namely 
$(\boldsymbol{V}, \boldsymbol{H}, \boldsymbol{V}^*)$ is a standard {H}ilbert triplet.

Under this setting, we define the solution of {\rm (P)} as follows.

\paragraph{Definition 2.1.} 
{\it The triplet $(\boldsymbol{u}, \boldsymbol{\mu }, \boldsymbol{\xi})$ 
is called the solution of {\rm (P)} if 
$\boldsymbol{u}=(u,u_\Gamma )$, $\boldsymbol{\mu }=(\mu,\mu _\Gamma)$, 
$\boldsymbol{\xi}=(\xi, \xi_\Gamma)$ satisfy
\begin{gather*}
	u \in H^1(0,T;H) \cap C \bigl( [0,T] ;V \bigr) \cap L^2(0,T;W), 
	\nonumber 
	\\ 
	u_\Gamma \in H^1(0,T;V_\Gamma^*) \cap L^\infty (0,T;V_\Gamma) \cap L^2(0,T;W_\Gamma), 
	\nonumber 
	\\
	\mu  \in L^2(0,T;V),
	\quad 
	\mu _\Gamma \in L^2(0,T;V_\Gamma ), 
	\\
	\xi \in L^2(0,T;H),
	\quad 
	\xi _\Gamma \in L^2(0,T;H_\Gamma ), 
	\nonumber 
	\\
	u_{|_\Sigma }=u_{\Gamma }, \quad 
	\mu _{|_\Sigma } = \mu _{\Gamma }
	\quad \mbox{a.e.\ on } \Sigma,
	\\
	\xi \in \beta (u) \quad \mbox{a.e.\ in }Q, \quad 
	\xi _\Gamma \in \beta _\Gamma (u_\Gamma)\quad \mbox{a.e.\ on }\Sigma
	\nonumber 
\end{gather*}
and solve
\begin{equation} 
	\bigl \langle u_\Gamma '(t), z_\Gamma
	\bigr \rangle _{V_\Gamma^*, V_\Gamma}
	+\int_{\Omega }^{} \nabla \mu(t)  \cdot \nabla z dx 
	+\int_{\Gamma }^{} \nabla _\Gamma \mu  _\Gamma(t) \cdot \nabla _\Gamma z_\Gamma d\Gamma 
	=0 
	\label{Def1} 
\end{equation}
for all $\boldsymbol{z}:=(z,z_\Gamma ) \in \boldsymbol{V}$, 
for a.a.\ $t\in (0,T)$, and
\begin{gather} 
	\tau \partial _t u-\Delta u + \xi + \pi (u) =f
	\quad \mbox{a.e.\ in }Q,
	\label{Def2}
	\\
	\mu _\Gamma =\partial_{\boldsymbol{\nu }} u - \Delta _\Gamma u_\Gamma 
	+\xi _\Gamma + \pi _\Gamma (u_\Gamma )-f_\Gamma
	\quad \mbox{a.e.\ on }\Sigma,
	\label{Def3}
	\\
	u(0)=u_0 
	\quad 
	\mbox{a.e.\ in }\Omega, \quad 
	u_\Gamma  (0)=u_{0\Gamma}
	\quad \mbox{a.e.\ on }\Gamma.
	\label{Def4}
\end{gather} }

\paragraph{Remark 1.} Taking $\boldsymbol{z}:=(1,1)$ in the weak formulation \eqref{Def1}, 
we see that \eqref{Def1} and \eqref{Def4} imply 
the mass conservation \eqref{mass} on the boundary. 
Moreover, for any $z \in {\mathcal D}(\Omega )$, 
taking $\boldsymbol{z}:=(z,0)$ in \eqref{Def1} and using 
the trace condition on $\mu $, we 
deduce that 
\begin{equation*}
	\begin{cases}
	-\Delta \mu(t) =0 \quad \mbox{a.e.\ in } \Omega, \\ 
	\mu _{|_\Gamma }(t)=\mu _\Gamma(t) \quad \mbox{a.e.\ on } \Gamma,
	\end{cases} 
	\quad \mbox{for a.a.\ } t \in (0,T),
\end{equation*}
whence 
the regularities
$\mu  \in L^2(0,T;V)$,
$\mu _\Gamma \in L^2(0,T;V_\Gamma )$ 
allow us to infer the higher regularity 
{\takeshi $\mu \in L^2(0,T;H^{3/2}(\Omega ))$. Moreover, 
the boundedness $\Delta \mu~(=0)$ in $L^2(0,T;H)$ gives us   
the property 
$\partial _{\boldsymbol{\nu }} \mu \in L^2(0,T;H_\Gamma)$}, 
as well as 
\begin{equation*} 
	\bigl \langle u_\Gamma '(t)+\partial _{\boldsymbol{\nu }} \mu(t) , z_\Gamma
	\bigr \rangle _{V_\Gamma^*, V_\Gamma}
	+ \int_{\Gamma}^{} \nabla _\Gamma \mu_\Gamma (t) \cdot \nabla _\Gamma z_\Gamma d\Gamma 
	=0 
	\quad \mbox{for all } z_\Gamma \in V_\Gamma,
\end{equation*} 
this is the weak formulation of \eqref{QS3}.

\subsection{Main theorems}

The first result states the existence of the solution. 
To the aim, we assume that: 
\begin{enumerate}
 \item[(A1)] $f \in L^2(0,T;H)$ and $f_\Gamma \in W^{1,1}(0,T;H_\Gamma)$; 
 \item[(A2)] $\boldsymbol{u}_0 := (u_0,u_{0\Gamma }) \in \boldsymbol{V}$;
 \item[(A3)] $\beta $, $\beta _{\Gamma }$, maximal monotone graphs in 
$\mathbb{R} \times \mathbb{R}$, that are the subdifferentials 
\begin{gather*}
	\beta =\partial \widehat{\beta}, \quad \beta _{\Gamma }
	=\partial \widehat{\beta }_{\Gamma }
\end{gather*}
of some proper lower semicontinuous and convex functions 
$\widehat{\beta }$ and $\widehat{\beta }_{\Gamma }: \mathbb{R} \to [0,+\infty ]$ 
satisfying $\widehat{\beta }(0)=\widehat{\beta}_{\Gamma }(0)=0$
with some 
effective domains $D(\widehat{\beta }) \supset D(\beta )$ and 
$D(\widehat{\beta }_\Gamma  ) \supset D(\beta _\Gamma)$, respectively. 
This implies that 
$0 \in \beta (0)$ and $0 \in \beta _{\Gamma }(0)$; 
 \item[(A4)]
$\pi $, $\pi _{\Gamma }: \mathbb{R} \to \mathbb{R}$ are {L}ipschitz continuous functions 
with {L}ipschitz constants $L$ and $L_{\Gamma}$, and they satisfy $\pi (0)=\pi _\Gamma (0)=0$; 
 \item[(A5)] $D(\beta _\Gamma ) \subseteq D(\beta )$ and 
there exist positive constants $\varrho, c_0 >0$ such that 
\begin{gather} 
	\bigl |\beta ^\circ (r) \bigr| 
	\le \varrho \bigl |\beta _{\Gamma }^\circ (r) \bigr |+c_0
	\quad 
	\mbox{for all } r \in D(\beta _\Gamma ),
	\label{A5}
\end{gather} 
where $\beta ^\circ $ and $\beta _\Gamma ^\circ $ denote the minimal sections of 
$\beta $ and $\beta _\Gamma $; 
  \item[(A6)] $m_\Gamma  \in {\rm int}D(\beta _\Gamma )$ and the compatibility conditions
$\widehat \beta (u_0) \in L^1(\Omega )$, 
$\widehat \beta _{\Gamma }(u_{0\Gamma}) \in L^1(\Gamma)$ hold.
\end{enumerate}
The minimal section $\beta ^\circ $ of $\beta $ is specified by
$\beta ^\circ (r):=\{ r^* \in \beta (r) : |r^*|=\min _{s \in \beta (r)} |s| \}$
and same definition applies to $\beta _{\Gamma }^\circ $. 
These assumptions are same as in \cite{CC13, CF15b}. 
{\colli Concerning assumption (A1), let us note that 
the regularity conditions for $f$ and $f_\Gamma $ are not symmetric: in fact, we need more 
regularity for the source term on the boundary since the equation on the boundary is of Cahn--Hilliard type, while the equation on the bulk turns out to be of the simpler Allen--Cahn type. 
Of course, here the condition $\tau>0$ plays a role, if the term $\tau \partial_t u$ would not be present in \eqref{QS2}, we would certainly need higher regularity for $f$.}

\paragraph{Theorem 2.1.} 
{\it Under the assumptions {\rm (A1)}--{\rm (A6)}, 
there exists a solution of the problem~{\rm (P)}.}

\bigskip
The second result states the continuous dependence on the data. 
The uniqueness of the component $\boldsymbol{u}$ of 
the solution is obtained from this theorem. 
Here, we {\colli just use the following regularity properties on the data:}
\begin{enumerate}
 \item[(A1)$^\prime $] $f \in L^2(0,T;V^*)$ and $f_\Gamma \in L^2(0,T;V_\Gamma ^*)$; 
 \item[(A2)$^\prime $] $u_0 \in H$ and $u_{0\Gamma} \in V_\Gamma^*$.
\end{enumerate}

Then, we obtain the continuous dependence on the data as follows:

\paragraph{Theorem 2.2.} 
{\it {\colli Under the assumption {\rm (A3)}--{\rm (A4)}, let, for $i=1,2$, 
$\boldsymbol{f}^{(i)}$, $\boldsymbol{u}_0^{(i)}$ satisfy {\rm (A1)$^\prime $}, {\rm (A2)$^\prime $} and assume that the corresponding solutions  $(\boldsymbol{u}^{(i)}, \boldsymbol{\mu }^{(i)}, \boldsymbol{\xi}^{(i)})$ exist.} Then, there 
{\colli is} a positive constant $C>0$, depending on $L$, $L_{\Gamma}$ and $T$, such that 
\begin{align}
	& 
	\bigl| u^{(1)}(t)-  u^{(2)}(t) 
	\bigr|_{H}^2 
	+
	\bigl| u_{\Gamma }^{(1)}(t)-  u_{\Gamma }^{(2)}(t) 
	\bigr|_{V_\Gamma ^*}^2 
	+ \int_{0}^{t} \bigl| u^{(1)}(s)-  u^{(2)}(s) 
	\bigr|_{V}^2 ds 
	 \nonumber \\
	& \quad {} + \int_{0}^{t} \bigl| u^{(1)}_\Gamma (s)-  u^{(2)}_\Gamma (s) 
	\bigr|_{V_\Gamma }^2 ds 
	\nonumber \\
	& \quad \quad \le C \left\{ 
	\bigl| u^{(1)}_0-  u^{(2)}_0 
	\bigr|_{H}^2 
	+
	\bigl| u_{0\Gamma }^{(1)}-  u_{0\Gamma }^{(2)}
	\bigr|_{V_\Gamma^*}^2 
	+
	\int_{0}^{t} 
	\bigl| f^{(1)}(s ) - f^{(2)}(s ) 
	\bigr|_{V^*}^2 ds 
	\right. \nonumber \\
	&  \left. 
	 \quad \quad \quad {} + \int_{0}^{t} 
	\bigl| f^{(1)}_\Gamma (s) - f^{(2)}_\Gamma (s) 
	\bigr|_{V_\Gamma^*}^2 ds 
	\right\}
	\label{conti}
\end{align} 
for all $t \in [0,T]$. }

\section{Approximate problem and uniform estimates}
\setcounter{equation}{0}

In this section, we first consider an approximate problem for {\rm (P)}, and then 
we obtain the uniform estimates. 
For each $\varepsilon, \lambda  \in (0,1]$, 
we introduce an approximate problem {\rm (P;$\lambda, \varepsilon$)} 
where the proof of the 
well-posedness of {\rm (P;$\lambda, \varepsilon$)} is given in Appendix.

\subsection{Moreau--Yosida regularization}

For each $\lambda  \in (0,1]$, we define 
$\beta _\lambda , \beta _{\Gamma,\lambda }:\mathbb{R} \to \mathbb{R}$, 
along with the associated resolvent operators 
$J_\lambda , J_{\Gamma,\lambda}:\mathbb{R} \to \mathbb{R}$
by 
\begin{gather*}
	\beta _\lambda (r)
	:= \frac{1}{\lambda } \bigl( r-J_\lambda (r) \bigr)
	:=\frac{1}{\lambda }\bigl( r-(I+\lambda  \beta )^{-1} (r)\bigr),
	\\
	\beta _{\Gamma, \lambda } (r)
	:= \frac{1}{\lambda  \varrho} \bigl( r-J_{\Gamma,\lambda }(r) \bigr )
	:=\frac{1}{\lambda  \varrho}\bigl( r- (I+\lambda \varrho \beta _\Gamma )^{-1} (r) \bigr)
\end{gather*}
for all $r \in \mathbb{R}$,
where $\varrho>0$ is the same constant as in the assumption \eqref{A5}. 
Note that the two definitions are not symmetric since in the second it is 
$\lambda \varrho$ and not directly $\lambda $ to be used as approximation parameter; 
let us note that 
this adaptation comes from the previous work \cite{CC13}. 
Now, we easily have 
$\beta _\lambda (0)=\beta _{\Gamma, \lambda }(0)=0$. 
Moreover, the related {M}oreau--{Y}osida regularizations $\widehat{\beta }_\lambda , 
\widehat{\beta }_{\Gamma,\lambda }$
of $\widehat{\beta }, \widehat{\beta}_{\Gamma }:\mathbb{R} \to \mathbb{R}$ fulfill
\begin{gather*}
	\widehat{\beta }_{\lambda }(r)
	:=\inf_{s \in \mathbb{R}}\left\{ \frac{1}{2\lambda } |r-s|^2
	+\widehat{\beta }(s) \right\} 
	= 
	\frac{1}{2\lambda } 
	\bigl| r-J_\lambda  (r) \bigr|^2+\widehat{\beta }\bigl (J_\lambda (r) \bigr )
	= \int_{0}^{r} \beta _\lambda (s)ds,
	\\
	\widehat{\beta }_{\Gamma, \lambda }(r)
	:=\inf_{s \in \mathbb{R}}\left\{ \frac{1}{2\lambda  \varrho } |r-s|^2
	+\widehat{\beta }_\Gamma (s) \right\} 
	= \int_{0}^{r} \beta _{\Gamma,\lambda } (s)ds
	\quad \mbox{for all } r\in \mathbb{R}.
\end{gather*}
It is well known that $\beta_\lambda $ is {L}ipschitz continuous with {L}ipschitz constant 
$1/\lambda $ 
and $\beta_{\Gamma, \lambda }$ is also {L}ipschitz continuous with constant 
$1/(\lambda  \varrho)$. In addition, for each $\lambda \in (0,1]$, we have the standard properties
\begin{gather}
	\bigl |\beta _\lambda (r) \bigr | \le \bigl |\beta ^\circ (r) \bigr |, \quad 
	\bigl |\beta _{\Gamma ,\lambda }(r) \bigr | \le \bigl |\beta _{\Gamma }^\circ (r) \bigr |,
	\nonumber \\
	0 \le \widehat{\beta }_\lambda (r) \le \widehat{\beta }(r), \quad 
	0 \le \widehat{\beta }_{\Gamma,\lambda } (r) \le \widehat{\beta }_{\Gamma }(r)
	\quad \mbox{for all } r \in \mathbb{R}.
	\label{prim}
\end{gather}
Let us point out that, using \cite[Lemma~4.4]{CC13}, we have 
\begin{gather}
	\bigl |\beta_\lambda (r)\bigr | 
	\le \varrho \bigl |\beta _{\Gamma,\lambda } (r)\bigr |+c_0
	\quad 
	\mbox{for all } r \in \mathbb{R},
	\label{A5e}
\end{gather} 
for all $\lambda \in (0,1]$
with the same constants $\varrho $ and $c_0$ as in \eqref{A5}. 

\medskip
Now for each $\varepsilon \in (0,1]$, 
letting $\boldsymbol{f}_\varepsilon :=(f_\varepsilon, f_{\Gamma, \varepsilon})$ 
and $\boldsymbol{u}_{0,\varepsilon }:=(u_{0,\varepsilon }, u_{0\Gamma, \varepsilon})$ be 
smooth approximations for 
$\boldsymbol{f}$ and $\boldsymbol{u}_0$, so that 
$\boldsymbol{f}_\varepsilon \in H^1(0,T;\boldsymbol{H})$ 
with $\boldsymbol{f}_\varepsilon (0) \in \boldsymbol{V}$ 
and $\boldsymbol{u}_{0,\varepsilon } \in \boldsymbol{W} \cap \boldsymbol{V} $ with 
$( -\Delta u_{0, \varepsilon }, 
\partial _{\boldsymbol{\nu }} u_{0, \varepsilon } - \Delta _\Gamma u_{0\Gamma, \varepsilon} ) 
\in \boldsymbol{V}$ 
satisfying 
\begin{gather}
	\boldsymbol{f}_\varepsilon \to \boldsymbol{f} 
	\quad \mbox{strongly in } L^2(0,T;\boldsymbol{H})
	\quad \mbox{as } \varepsilon \to 0,
	\label{apf1} \\
	\quad | f_{\Gamma, \varepsilon } - f_\Gamma |_{L^2(0,T;H_\Gamma )} \le \varepsilon ^{1/2} C_0 
	\quad \mbox{for all } \varepsilon \in (0,1], 
	\label{apf2} 
	\\
	\boldsymbol{u}_{0,\varepsilon } \to \boldsymbol{u}_0
	\quad \mbox{strongly in } \boldsymbol{V}
	\quad \mbox{as } \varepsilon \to 0,
	\label{apu}
	\\
	\int_{\Omega }^{} \widehat{\beta }_\lambda (u_{0, \varepsilon }) dx 
	\le C_0, \quad 
	\int_{\Gamma }^{} \widehat{\beta }_{\Gamma, \lambda }(u_{0\Gamma, \varepsilon }) d\Gamma  \le 
	\left( 1+\frac{\varepsilon ^{1/2}}{\lambda }\right) C_0
	\quad \mbox{for all } \varepsilon \in (0,1],
	\label{apbound} 
\end{gather}
where $C_0$ is a positive constant independent of $\varepsilon, \lambda \in (0,1]$. 
Indeed, $\boldsymbol{f}_\varepsilon $ and $\boldsymbol{u}_{0, \varepsilon} $ 
satisfying \eqref{apf1}--\eqref{apbound} are given in the Appendix.  
Hereafter we use $C^*:=(1+\varepsilon ^{1/2}/\lambda )^{1/2}$ which 
{\revis satisfies} $C^* \searrow 1$ as $\varepsilon \searrow 0$ for $\lambda \in (0,1]$. 
Then we can solve 
the following auxiliary problem.

\paragraph{Proposition 3.1.} 
{\it Under the assumptions {\rm (A1)}--{\rm (A6)}, 
for each $\varepsilon, \lambda  \in (0,1]$ there 
exists a unique pair 
\begin{gather}
	\boldsymbol{u}_{\lambda, \varepsilon} 
	:=(u_{\lambda, \varepsilon}, u_{\Gamma, \lambda, \varepsilon}) 
	\in W^{1,\infty }(0,T;\boldsymbol{H}) 
	\cap H^1(0,T;\boldsymbol{V}) 
	\cap L^\infty (0,T;\boldsymbol{W}), 
	\label{pro311} \\
	\boldsymbol{\mu }_{\lambda, \varepsilon} :=(\mu_{\lambda, \varepsilon}, \mu _{\Gamma, \lambda, \varepsilon} ) 
	\in L^\infty (0,T;\boldsymbol{V}) \cap L^2(0,T;\boldsymbol{W})
	\label{pro312}
\end{gather}
satisfying
\begin{gather} 
	\varepsilon \partial _t u_{\lambda, \varepsilon}
	-\Delta \mu _{\lambda, \varepsilon} = 0
	\quad \mbox{a.e.\ in }Q,
	\label{apDef1} 
	\\
	\varepsilon \mu _{\lambda, \varepsilon} 
	=\tau \partial _t u_{\lambda, \varepsilon}
	-\Delta u_{\lambda, \varepsilon} 
	+ \beta _\lambda (u_{\lambda, \varepsilon}) 
	+ \pi (u_{\lambda, \varepsilon}) -f_\varepsilon 
	\quad \mbox{a.e.\ in }Q,
	\label{apDef2}
	\\
	\partial_t u_{\Gamma, \lambda, \varepsilon}
	+ \partial _{\boldsymbol{\nu }} \mu _{\lambda, \varepsilon}
	- \Delta _\Gamma \mu _{\Gamma, \lambda, \varepsilon} = 0
	\quad \quad \mbox{a.e.\ on }\Sigma,
	\label{apDef3}
	\\
	\mu _{\Gamma, \lambda, \varepsilon}=
	\varepsilon \partial _t u_{\Gamma, \lambda, \varepsilon}
	+ \partial_{\boldsymbol{\nu }} u_{\lambda, \varepsilon} - \Delta _\Gamma u_{\Gamma, \lambda, \varepsilon}
	+ \beta _{\Gamma, \lambda}(u_{\Gamma, \lambda, \varepsilon}) 
	+ \pi _\Gamma (u_{\Gamma, \lambda, \varepsilon})
	- f_{\Gamma, \varepsilon }
	\quad \mbox{a.e.\ on }\Sigma,
	\label{apDef4}
	\\
	u_{\lambda, \varepsilon}(0)=u_{0, \varepsilon } 
	\quad 
	\mbox{a.e.\ in }\Omega, \quad 
	u_{\Gamma, \lambda, \varepsilon} (0)=u_{0\Gamma, \varepsilon }
	\quad \mbox{a.e.\ on }\Gamma.
	\label{apDef5}
\end{gather}} 
\indent
As a remark, the trace conditions is included in the 
regularities \eqref{pro311}--\eqref{pro312}. The strategy of the proof of this proposition 
is based on the previous work \cite[Theorems~2.2, 4.2]{CF15b}. 
It is given in the Appendix.

\subsection{A priori estimates}

In order to obtain the uniform estimates independent of the 
approximate parameter $\varepsilon $ and $\lambda $, we use the following 
type of Poincar\'e--Wirtinger inequalities (see, e.g., \cite{Nec67, Heb96}): 
there exists a positive constant $C_P$ such that 
\begin{gather} 
	|z|_H^2 
	\le C_P \left\{ 
	\int_{\Omega }^{} |\nabla z| ^2 dx + \int_{\Gamma }^{} |z_{|_\Gamma }|^2 d\Gamma 
	\right\}  
	\quad \mbox{for all } z \in V,
	\label{poin1}
	\\
	|z_\Gamma |_{H_\Gamma }^2 
	\le C_P 
	\int_{\Gamma }^{} |\nabla _\Gamma z_\Gamma |^2d\Gamma 
	\quad 
	\mbox{for all }z_\Gamma \in V_\Gamma \ \mbox{with } \int_{\Gamma }^{} z_\Gamma d\Gamma =0,
	\label{poin2}
	\\ 
	|\boldsymbol{z}|_{\boldsymbol{V}} ^2 
	\le C_P \left\{ 
	\int_{\Omega }^{}|\nabla z|^2 dx + \int_{\Gamma }^{}|\nabla _{\Gamma } z_\Gamma |^2d\Gamma 
	\right\} 
	\quad \mbox{for all } \boldsymbol{z} \in \boldsymbol{V} \ \mbox{with } \int_{\Gamma }^{} z_\Gamma d\Gamma =0.
	\label{poin3}
\end{gather}
Moreover, from the compactness inequality recalled in 
\cite[Chapter~1, Lemme 5.1]{Lio69} 
or \cite[Section~8, Lemma~8]{Sim87}, 
for each $\delta >0$ there exists 
a positive constant $C_\delta$ depending on $\delta $ such that
\begin{gather} 
	{\fukao |z|_{H^{1-\alpha }(\Omega)}^2 \le \delta |z|_{V}^2 + C_\delta |z|_{H}^2 
	\quad \mbox{for all } z \in V \ \mbox{and~all~} \alpha  \in (0,1),}
	\label{cp} \\
	|z_\Gamma |_{H_\Gamma }^2 \le \delta |z_\Gamma |_{V_\Gamma}^2 + C_\delta |z_\Gamma |_{V_\Gamma ^*}^2 
	\quad \mbox{for all } z_\Gamma \in V_\Gamma,
	\label{cp0}
\end{gather}
because we have the compact embeddings 
{\fukao $V \mathop{\hookrightarrow} \mathop{\hookrightarrow} H^{1-\alpha }(\Omega )
\mathop{\hookrightarrow} \mathop{\hookrightarrow} H$ for $\alpha \in (0,1)$}, 
(see, e.g., \cite[Chapter~1, Theorem~16.1]{LM70}) and 
$V_\Gamma \mathop{\hookrightarrow} \mathop{\hookrightarrow} 
H_\Gamma \mathop{\hookrightarrow} \mathop{\hookrightarrow} V_\Gamma ^*$, respectively.
Next, from the standard theorem for the 
trace {\fukao operators $\gamma_0 (z):= z_{|_\Gamma}$ of 
$\gamma_0:V \to H^{1/2}(\Gamma )$ and 
$\gamma_0:H^{1-\alpha }(\Omega ) \to H^{(1-\alpha )-1/2}(\Gamma) \subset H_\Gamma$ for $\alpha  \in (0,1/2)$,}
(see, e.g.,  \cite[Chapter~2, Theorem~2.24]{BG87}, \cite[Chapter~2, Theorem~5.5]{Nec67}), 
we see that 
there exists a positive constant $c_1$ such that 
\begin{gather} 
	|z_\Gamma |_{H_\Gamma }^2 \le |z_\Gamma |_{H^{1/2}(\Gamma )}^2 
	\le c_1 |z|_{V}^2,
	\quad {\fukao 
	|z_\Gamma |_{H_\Gamma }^2 \le c_1 |z|_{H^{1-\alpha }(\Omega )}^2
	\quad \mbox{for all } \boldsymbol{z} \in \boldsymbol{V},}
	\label{co}
\end{gather}
{\fukao for $\alpha  \in (0,1/2)$}.
Moreover, the boundedness of some recovering operator ${\mathcal R}:H^{1/2}(\Gamma ) \to V$ 
of the trace $\gamma_0 $
gives us 
\begin{equation}
	\bigl| {\mathcal R} z_\Gamma  \bigr|_{V}^2  
	\le c_2 |z_\Gamma |_{H^{1/2}(\Gamma )}^2
	\le c_2 |z_\Gamma |_{V_\Gamma }^2
	\quad \mbox{for all } z_\Gamma \in V_\Gamma,
	\label{co0}
\end{equation}
where $c_2$ is a positive constant (see, e.g., \cite[Chapter~2, Theorem~2.24]{BG87}, \cite[Chapter~2, Theorem~5.7]{Nec67}).

\paragraph{Lemma 3.1.}
{\it There exist two positive constants $M_1$ and 
$M_2$, 
depending on $\tau $ but independent of $\varepsilon $ and $\lambda  \in (0,1]$, 
such that}
\begin{gather} 
	|u_{\lambda, \varepsilon }|_{L^\infty (0,T;V)} \le C^* M_1,
	\quad 
	|u_{\Gamma, \lambda, \varepsilon }|_{L^\infty (0,T;H_\Gamma )} \le c_1^{1/2} C^* M_1,
	\label{m1}\\
	|\partial _t u_{\lambda, \varepsilon }|_{L^2(0,T;H)}
	+
	\varepsilon ^{1/2} |\partial _t u_{\Gamma, \lambda, \varepsilon }|_{L^2(0,T;H_\Gamma)}
	+{\revis 
	|\nabla _\Gamma u_{\Gamma, \lambda,\varepsilon}|_{L^\infty (0,T;H_\Gamma)}}
	+
	\bigl|
	\widehat{\beta }_{\lambda } ( u_{\lambda, \varepsilon}) 
	\bigr|_{L^\infty (0,T;L^1(\Omega ))}
	\nonumber \\
	\quad {}
	+
	\bigl|
	\widehat{\beta }_{\Gamma, \lambda } ( u_{\Gamma, \lambda, \varepsilon} ) 
	\bigr|_{L^\infty (0,T;L^1(\Gamma))}
	+{\revis 
	|\nabla \mu _{\lambda, \varepsilon}|_{L^2(0,T;H)}}
	+{\revis 
	|\nabla_\Gamma  \mu _{\Gamma, \lambda, \varepsilon}|_{L^2(0,T;H_{\Gamma })}}
	\le C^* M_2.
	\label{m2}
\end{gather}

\paragraph{Proof.} We test \eqref{apDef2} at time $s$ by 
$\partial_s u_{\lambda, \varepsilon}$, 
the time derivative of $u_{\lambda, \varepsilon}$. Integrating the resultant over $\Omega \times (0,t)$ 
leads to
\begin{align}
	& \tau \int_{0}^{t} 
	\bigl| \partial_s u_{\lambda, \varepsilon}(s) \bigr|^2_H ds
	+ {\revis 
	\frac{1}{2}
	\bigl| \nabla u_{\lambda, \varepsilon}(t) \bigr|_{H}^2}
	+
	\int_{\Omega }^{}
	\widehat{\beta }_{\lambda } \bigl( u_{\lambda, \varepsilon}(t) \bigr) 
	dx
	+
	\int_{\Omega }^{}
	\widehat{\pi } \bigl( u_{\lambda, \varepsilon}(t) \bigr) dx 
	\nonumber \\
	& \quad 
	{} - 
	\int_{0}^{t} \bigl( \partial _{\boldsymbol{\nu }} u_{\lambda, \varepsilon}(s), 
	\partial_s u_{\Gamma, \lambda, \varepsilon}(s) 
	\bigr)_{H_\Gamma } ds
	 - 
	 \int_{0}^{t} 
	 \bigl( f_\varepsilon (s), \partial_s u_{\lambda, \varepsilon}(s) \bigr)_{H} ds 
	\nonumber \\
	& = 
	\varepsilon \int_{0}^{t} 
	\bigl(
	\mu _{\lambda,\varepsilon  }(s), 
	 \partial_s u_{\lambda, \varepsilon}(s) \bigr)_{H} ds 
	 + {\revis 
	\frac{1}{2}
	| \nabla u_{0, \varepsilon } |_{H}^2}
	+
	\int_{\Omega }^{}
	\widehat{\beta }_{\lambda } ( u_{0, \varepsilon } ) 
	dx
	+
	\int_{\Omega}^{}
	\widehat{\pi } ( u_{0, \varepsilon } ) 
	dx 
	\label{l3.1.a}
\end{align} 
for all $t \in [0,T]$. 
Next, testing \eqref{apDef4} by 
$\partial_s u_{\Gamma, \lambda, \varepsilon}$ and 
integrating the resultant over $\Gamma \times (0,t)$, we obtain 
\begin{align}
	& -\int_{0}^{t} \bigl( \partial _{\boldsymbol{\nu }} u_{\lambda, \varepsilon}(s), 
	\partial_s u_{\Gamma, \lambda, \varepsilon}(s) 
	\bigr)_{H_\Gamma } ds 
	\nonumber \\
	& = \varepsilon \int_{0}^{t} 
	\bigl|
	\partial_s u _{\Gamma, \lambda, \varepsilon}(s)
	\bigr|_{H_\Gamma }^2 ds 
	- \int_{0}^{t} \bigl( \mu _{\Gamma, \lambda, \varepsilon}(s), 
	\partial_s u_{\Gamma, \lambda, \varepsilon}(s) \bigr)_{H_\Gamma } ds 
	+ {\revis \frac{1}{2} \bigl| \nabla u_{\Gamma, \lambda, \varepsilon}(t) \bigr|_{H_\Gamma}^2}
	- {\revis \frac{1}{2} | \nabla u_{0\Gamma, \varepsilon }|_{H_\Gamma}^2}
	\nonumber \\
	& \quad {}
	+
	\int_{\Gamma}^{}
	\widehat{\beta }_{\Gamma, \lambda } \bigl( u_{\Gamma, \lambda, \varepsilon}(t) \bigr) 
	d\Gamma 
	+
	\int_{\Gamma}^{}
	\widehat{\pi }_{\Gamma} \bigl( u_{\Gamma, \lambda, \varepsilon}(t) \bigr) 
	d\Gamma 
	\nonumber \\
	& \quad {}
	- 
	\int_{\Gamma }^{}
	\widehat{\beta }_{\Gamma, \lambda } ( u_{0\Gamma, \varepsilon } ) 
	d\Gamma
	-
	\int_{\Gamma }^{}
	\widehat{\pi }_{\Gamma } ( u_{0\Gamma, \varepsilon } ) 
	d\Gamma 
	- \int_{0}^{t} \bigl( f_{\Gamma, \varepsilon}(s), \partial_s u_{\Gamma, \lambda, \varepsilon}(s) \bigr)_{H_\Gamma } ds 
	\label{l3.1.b}
\end{align}
for all $t \in [0,T]$. 
On the other hand, 
testing \eqref{apDef1} by 
$\mu _{\lambda, \varepsilon}$, testing \eqref{apDef3} by 
$\mu _{\Gamma, \lambda, \varepsilon}$, and adding them, we infer that 
\begin{align}
	& \varepsilon \int_{0}^{t} 
	\bigl(
	\mu _{\lambda,\varepsilon  }(s), 
	 \partial_s u_{\lambda, \varepsilon}(s) \bigr)_{H} ds 
	 + 
	\int_{0}^{t} \bigl( \mu _{\Gamma, \lambda, \varepsilon}(s), 
	\partial_s u_{\Gamma, \lambda, \varepsilon}(s) \bigr)_{H_\Gamma } ds 
	\nonumber \\
	& = - {\revis 
	\int_{0}^{t} 
	\bigl| \nabla \mu _{\lambda, \varepsilon} (s) \bigr|_{H}^2 ds 
	}
	- {\revis 
	 \int_{0}^{t} \bigl| \nabla_{\Gamma}  \mu _{\Gamma, \lambda, \varepsilon} (s) \bigr|_{H_{\Gamma}}^2 ds }
	 \label{l3.1.c}
\end{align} 
for all $t \in [0,T]$. 
Combining \eqref{l3.1.a}--\eqref{l3.1.c} and using \eqref{prim}, we have 
\begin{align}
	& \tau \int_{0}^{t} 
	\bigl| \partial_s u_{\lambda, \varepsilon}(s) \bigr|^2_H ds
	+
	\varepsilon \int_{0}^{t} 
	\bigl|
	\partial_s u _{\Gamma, \lambda, \varepsilon}(s)
	\bigr|_{H_\Gamma }^2 ds 
	+ {\revis 
	\frac{1}{2}
	\bigl| \nabla u_{\lambda, \varepsilon}(t) \bigr|_{H}^2 }
	+
	\int_{\Omega }^{}
	\widehat{\beta }_{\lambda } \bigl( u_{\lambda, \varepsilon}(t) \bigr) 
	dx
	\nonumber \\
	& \quad 
	{} 
	+ 
	{\revis \frac{1}{2} \bigl| \nabla u_{\Gamma, \lambda, \varepsilon}(t) \bigr|_{H_\Gamma}^2
	}
	+ \int_{\Gamma }^{}
	\widehat{\beta }_{\Gamma, \lambda } \bigl( u_{\Gamma, \lambda, \varepsilon}(t) \bigr) 
	d\Gamma
	+ {\revis 
	\int_{0}^{t} 
	\bigl| \nabla \mu _{\lambda, \varepsilon} (s) \bigr|_{H}^2 ds 
	}
	+ {\revis 
	 \int_{0}^{t} \bigl| \nabla_{\Gamma}  \mu _{\Gamma, \lambda, \varepsilon} (s) \bigr|_{H_{\Gamma}}^2 ds 
	 }
	\nonumber \\
	& \le {\revis 
	\frac{1}{2}
	| \nabla u_{0, \varepsilon } |_{H}^2
	}
	+
	\int_{\Omega }^{}
	\widehat{\beta }_\lambda  ( u_{0,\varepsilon } ) 
	dx
	+ 
	\int_{\Omega }^{} 
	\bigl| \widehat{\pi } \bigl( u_{\lambda, \varepsilon}(t) \bigr) \bigr|
	dx 
	+
	\int_{\Omega }^{} \bigl| 
	\widehat{\pi } ( u_{0,\varepsilon} ) 
	\bigr|
	dx 
	 \nonumber \\
	& \quad {}
	+ {\revis 
	\frac{1}{2} | \nabla u_{0\Gamma, \varepsilon }|_{H_\Gamma}^2
	}
	+
	\int_{\Gamma }^{}
	\widehat{\beta }_{\Gamma, \lambda } ( u_{0\Gamma, \varepsilon } ) 
	d\Gamma
	+
	\int_{\Gamma }^{}
	\bigl| \widehat{\pi }_{\Gamma} \bigl( u_{\Gamma, \lambda, \varepsilon}(t) \bigr) \bigr|
	d\Gamma 
	+
	\int_{\Gamma }^{}
	\bigl| \widehat{\pi }_{\Gamma } ( u_{0\Gamma, \varepsilon } ) 
	\bigr| 
	d\Gamma 
	\nonumber \\
	& \quad {}
	{\colli {}+
	\int_{0}^{t} 
	 \bigl( f_\varepsilon (s), \partial_s u_{\lambda, \varepsilon}(s) \bigr)_{H} ds} 
	+ 
	\int_{0}^{t} \bigl( f_{\Gamma, \varepsilon }(s), \partial_s u_{\Gamma, \lambda, \varepsilon}(s) 
	\bigr)_{H_\Gamma } ds
	\label{l3.1.d}
\end{align} 
for all $t \in [0,T]$. 
Therefore, in order to estimate the right hand side of \eqref{l3.1.d}
we prepare the estimate of $|u_{\lambda, \varepsilon  }(s)|_{H}$. Indeed, 
from the {Y}oung inequality, 
we see that 
\begin{align}
	\frac{1}{2}
	\bigl|
	u_{\lambda, \varepsilon  }(s) 
	\bigr|_{H}^2 
	& = 
	\int_{0}^{t} 
	\bigl( \partial _s u_{\lambda, \varepsilon  }(s),
	u_{\lambda, \varepsilon}(s)  \bigr)_{H} ds 
	+
	\frac{1}{2}
	|
	u_{0, \varepsilon }
	|_{H}^2 
	\nonumber \\
	& \le \frac{\tilde{\delta }}{2} \int_{0}^{t} \bigl| \partial _s u_{\lambda, \varepsilon  }(s) 
	\bigr|_{H}^2 ds 
	+ \frac{1}{2\tilde{\delta }} \int_{0}^{t} \bigl| u_{\lambda, \varepsilon  }(s) 
	\bigr|_{H}^2 ds 
	+\frac{1}{2}
	|
	u_{0, \varepsilon }
	|_{H}^2 
	 \label{l3.1.e}
\end{align} 
for all $t \in [0,T]$ with some $\tilde{\delta }>0$. 
Now, from (A4) we can use 
the fact that $|\pi (r)| = | \pi (r)-\pi(0)| \le L |r|$, and then we deduce
\begin{equation*}
	\bigl| 
	\widehat{\pi }(r)
	\bigr| \le   
	\int_{0}^{r} \bigl| \pi(l ) \bigr| dl 
	\le \frac{L}{2}r^2
	\quad \mbox{for all } r \in \mathbb{R}.
\end{equation*}
Therefore, by taking $\tilde{\delta }:=\tau /(5L)$ in \eqref{l3.1.e} 
\begin{align}
	\int_{\Omega }^{} 
	\bigl| \widehat{\pi } \bigl( u_{\lambda, \varepsilon}(t) \bigr) \bigr|
	dx 
	& \le \frac{L}{2} \int_{\Omega }^{} 
	\bigl| u_{\lambda, \varepsilon}(t) \bigr|^2 
	dx 
	\nonumber \\
	& \le \frac{\tau}{10} \int_{0}^{t} \bigl| \partial _s u_{\lambda, \varepsilon  }(s) 
	\bigr|_{H}^2 ds 
	+  \frac{5L^2}{2\tau } \int_{0}^{t} \bigl| u_{\lambda, \varepsilon  }(s) 
	\bigr|_{H}^2 ds 
	+\frac{L}{2}
	|
	u_{0, \varepsilon }
	|_{H}^2 
	\label{cor1}
\end{align}
and analogously  
\begin{equation}
	\int_{\Omega }^{} \bigl| 
	\widehat{\pi } ( u_{0,\varepsilon } ) 
	\bigr|
	dx 
	\le \frac{L}{2} \int_{\Omega }^{} 
	| u_{0,\varepsilon } |^2 
	dx, \quad 
	\int_{\Gamma }^{}
	\bigl| \widehat{\pi }_{\Gamma } ( u_{0\Gamma,\varepsilon } ) 
	\bigr| 
	d\Gamma 
	\le 
	\frac{L_\Gamma }{2} \int_{\Gamma }^{}
	| u_{0\Gamma,\varepsilon } |^2
	d\Gamma.
	\label{cor2}
\end{equation}
Additionally, by using \eqref{cp}, \eqref{co} and \eqref{l3.1.e} with 
$\tilde{\delta }:=\tau /(10C_\delta )$ 
\begin{align}
	& \int_{\Gamma }^{}
	\bigl| \widehat{\pi }_{\Gamma} \bigl( u_{\Gamma, \lambda, \varepsilon}(t) \bigr) \bigr|
	d\Gamma 
	\nonumber \\
	& \quad \le 
	\frac{L_\Gamma }{2} \int_{\Gamma }^{}
	\bigl|  u_{\Gamma, \lambda, \varepsilon}(t) \bigr|^2
	d\Gamma 
	\nonumber \\
	& {\fukao \quad {} \le 
	\frac{c_1 L_\Gamma }{2} \bigl| u_{\lambda, \varepsilon}(t) \bigr|_{H^{1-\alpha }(\Omega)}^2} 
	\nonumber \\
	& \quad \le \delta \bigl| u_{\lambda, \varepsilon}(t) \bigr|_{V}^2 + 
	C_\delta \bigl| u_{\lambda, \varepsilon}(t) \bigr|_{H}^2
	\nonumber \\
	& \quad \le 
	\delta \bigl| u_{\lambda, \varepsilon}(t) \bigr|_{V}^2 + 
	\frac{\tau }{10}  \int_{0}^{t} \bigl| \partial _s u_{\lambda, \varepsilon  }(s) 
	\bigr|_{H}^2 ds 
	+ \frac{10 C_\delta^2 }{\tau } \int_{0}^{t} \bigl| u_{\lambda, \varepsilon  }(s) 
	\bigr|_{H}^2 ds 
	+C_\delta 
	|
	u_{0,\varepsilon }
	|_{H}^2 
	\label{cor3}
\end{align}
for some {\fukao $\alpha  \in (0,1/2)$ and} $\delta >0$. 
Moreover, from the {Y}oung inequality
\begin{align}
	\int_{0}^{t} 
	 \bigl( f_\varepsilon (s), \partial_s u_{\lambda, \varepsilon}(s) \bigr)_{H} ds 
	& \le 
	 \frac{\tau }{10} \int_{0}^{t} \bigl| \partial _s u_{\lambda, \varepsilon  }(s) 
	\bigr|_{H}^2 ds 
	+ \frac{5}{2\tau } \int_{0}^{t} \bigl| f_\varepsilon (s) 
	\bigr|_{H}^2 ds
	\label{cor4}, \\
	\int_{0}^{t} \bigl( f_{\Gamma,\varepsilon  }(s), \partial_s u_{\Gamma, \lambda, \varepsilon}(s) 
	\bigr)_{H_\Gamma } ds 
	& \le 
	\frac{\varepsilon  }{2} \int_{0}^{t} \bigl| \partial _s u_{\Gamma, \lambda, \varepsilon  }(s) 
	\bigr|_{H_\Gamma }^2 ds 
	+ \frac{1}{2\varepsilon } \int_{0}^{t} \bigl| f_{\Gamma, \varepsilon} (s) 
	\bigr|_{H_\Gamma }^2 ds
	\label{cor5}
\end{align} 
for all $t \in [0,T]$. 
Therefore, collecting \eqref{cor1}--\eqref{cor5}, 
adding $(1/2)|u_{\lambda,\varepsilon}(s)|_H^2$ to 
the both side of \eqref{l3.1.d} and using \eqref{l3.1.e} with $\tilde{\delta }:=\tau /5$, 
we obtain
\begin{align}
	& \frac{\tau }{2} \int_{0}^{t} 
	\bigl| \partial_s u_{\lambda, \varepsilon}(s) \bigr|^2_H ds
	+ \frac{\varepsilon }{2}\int_{0}^{t} 
	\bigl| \partial_s u_{\Gamma, \lambda, \varepsilon}(s) \bigr|^2_{H_\Gamma } ds
	+
	\frac{1}{2}
	\bigl| u_{\lambda, \varepsilon}(t) \bigr|_{V}^2
	+
	\int_{\Omega }^{}
	\widehat{\beta }_{\lambda } \bigl( u_{\lambda, \varepsilon}(t) \bigr) 
	dx
	\nonumber \\
	& \quad 
	{} 
	+ {\revis 
	\frac{1}{2} \bigl| \nabla u_{\Gamma, \lambda, \varepsilon}(t) \bigr|_{H_\Gamma}^2
	}
	+ 
	\int_{\Gamma }^{}
	\widehat{\beta }_{\Gamma, \lambda } \bigl( u_{\Gamma, \lambda, \varepsilon}(t) \bigr) 
	d\Gamma 
	+ {\revis 
	\int_{0}^{t} 
	\bigl| \nabla \mu _{\lambda, \varepsilon} (s) \bigr|_{H}^2 ds 
	}
	+
	{\revis 
	 \int_{0}^{t} \bigl| \nabla_{\Gamma}  \mu _{\Gamma, \lambda, \varepsilon} (s) \bigr|_{H_{\Gamma}}^2 ds 
	 }
	\nonumber \\
	& \le 
	\frac{5}{2\tau }\int_{0}^{t} \bigl| u_{\lambda, \varepsilon} (s) \bigr|_H^2 ds +
	\frac{1}{2}
	|u_{0, \varepsilon } |_{V}^2
	+
	\int_{\Omega }^{}
	\widehat{\beta }_\lambda ( u_{0, \varepsilon } ) 
	dx
	+ 
	\frac{5L^2}{2\tau} \int_{0}^{t} \bigl| u_{\lambda, \varepsilon  }(s) 
	\bigr|_{H}^2 ds 
	+L
	|
	u_{0,\varepsilon }
	|_{H}^2 
	\nonumber \\
	& \quad {}
	+{\revis 
	\frac{5}{2\tau } \int_{0}^{t} \bigl| f_\varepsilon (s) \bigr|_{H}^2 ds }
	+ \frac{1}{2} | u_{0\Gamma, \varepsilon }|_{V_\Gamma}^2
	+
	\int_{\Gamma }^{}
	\widehat{\beta }_{\Gamma,\lambda } ( u_{0\Gamma,\varepsilon } ) 
	d\Gamma
	+
	\delta \bigl| u_{\lambda, \varepsilon}(t) \bigr|_{V}^2 
	\nonumber \\
	& \quad {} 
	+ \frac{10C_\delta^2}{\tau } \int_{0}^{t} \bigl| u_{\lambda, \varepsilon} (s) \bigr|_H^2 ds 
	+ C_\delta |u_{0,\varepsilon }|_H^2
	+
	\frac{L_\Gamma }{2} |u_{0\Gamma,\varepsilon  }|_{H_\Gamma }^2
	+\frac{1}{2\varepsilon } \int_{0}^{t} \bigl| f_{\Gamma, \varepsilon} (s) 
	\bigr|_{H_\Gamma }^2 ds
	\label{l3.1.f}
\end{align} 
for all $t \in [0,T]$. 
Thus, taking $\delta :=1/4$, and using \eqref{apf1}--\eqref{apbound} and the {G}ronwall inequality,
we see that there exists a positive constant $M_1$
depending on 
$|u_{0} |_{V}$, $C_0$, 
$L$, 
$|f|_{L^2(0,T;H)}$, 
$|u_{0\Gamma}|_{V_\Gamma}$, 
$L_\Gamma$, 
$|f_\Gamma |_{L^2(0,T;H_\Gamma )}$
and $\tau $, in which $M_1 \to +\infty $ as 
$\tau \to 0$, 
independent of $\varepsilon $ and 
$\lambda $ such that 
\begin{equation*} 
	|u_{\lambda, \varepsilon }|_{L^\infty (0,T;V)} \le C^* M_1.
\end{equation*} 
Moreover, using this, \eqref{co} implies 
\begin{equation*} 
	|u_{\Gamma, \lambda, \varepsilon }|_{L^\infty (0,T;H_\Gamma )} \le c_1^{1/2} C^* M_1,
\end{equation*} 
and from \eqref{l3.1.f} we obtain the second estimates \eqref{m2}. 
\hfill $\Box$ 

\bigskip
By comparison, {\colli the following estimates for 
$\Delta \mu _{\lambda, \varepsilon}$ and $\boldsymbol{u}_{\lambda, \varepsilon}':=
(\partial _t u_{\lambda, \varepsilon}, \partial _t u_{\Gamma, \lambda, \varepsilon })$
are obtained.}

\paragraph{Lemma 3.2.}  
{\it Let 
$M_2$ be the same constant in Lemma~3.1. 
Then, for each $\varepsilon $ and $\lambda \in (0,1]$ 
the estimates
\begin{gather}
	|\Delta \mu _{\lambda, \varepsilon  }|_{L^2(0,T;H)} 
	\le \varepsilon C^* M_2,
	\label{m2b}
	\\
	| \boldsymbol{u}_{\lambda, \varepsilon}'
	|_{L^2(0,T;\boldsymbol{V}^*)} 
	\le 2 C^* M_2
	\label{m2c}
\end{gather}
hold.}
\paragraph{Proof.} 
From \eqref{apDef1}, we easily see that 
\begin{align*}
	|\Delta \mu _{\lambda, \varepsilon  }|_{L^2(0,T;H)} 
	& = {\revis |\varepsilon \partial _t u_{\lambda, \varepsilon }|_{L^2(0,T;H)} }
	\nonumber \\
	& \le \varepsilon C^* M_2.
\end{align*}
Next, we separate $\boldsymbol{u}_{\lambda, \varepsilon}'$ as follows:
\begin{equation*}
	\boldsymbol{u}_{\lambda, \varepsilon}' 
	= \bigl( (1-\varepsilon )\partial_t u_{\lambda, \varepsilon}, 0 \bigr)
	+
	( \varepsilon \partial_t u_{\lambda, \varepsilon}, \partial _t u_{\Gamma, \lambda, \varepsilon} ).
\end{equation*} 
Then, the first term of the right hand side is estimated as follows:
\begin{align}
	\bigl| \bigl( (1-\varepsilon )\partial_t u_{\lambda, \varepsilon}, 0 \bigr) \bigr|_{L^2(0,T;\boldsymbol{V}^*)} 
	& = \bigl| (1-\varepsilon ) \partial_t u_{\lambda, \varepsilon} \bigr|_{L^2(0,T;H)} 
	\nonumber \\
	& \le (1-\varepsilon ) C^* M_2
	\nonumber \\
	& \le C^*M_2. 
	\label{m2d}
\end{align}
Moreover, for the second term, we see that
\begin{equation}
	\bigl| 
	(\varepsilon \partial _t u_{\lambda, \varepsilon}, \partial _t u_{\Gamma, \lambda, \varepsilon}) 
	\bigr|_{L^2(0,T;\boldsymbol{V}^*)} \le C^*M_2. 
	\label{m2e}
\end{equation}  
Indeed, from \eqref{apDef1} and \eqref{apDef3}, we have 
\begin{gather*} 
	\int_{0}^{T} 
	\bigl( \varepsilon \partial _t u_{\lambda, \varepsilon} (t), \eta(t)  \bigr)_H
	dt
	+ \int_{0}^{T}
	\bigl( \partial _t u_{\Gamma, \lambda, \varepsilon} (t), \eta _\Gamma(t) \bigr)_{H_\Gamma}
	dt
	\\ 
	= {\revis \int_{0}^{T} \bigl( \nabla \mu_{\lambda, \varepsilon}(t), 
	\nabla \eta (t) \bigr)_{H} dt }
	+ {\revis
	\int_{0}^{T} \bigl( \nabla _\Gamma \mu _{\Gamma, \lambda, \varepsilon  }(t),
	\nabla _\Gamma \eta _\Gamma (t) \bigr) _{H_\Gamma} dt }
\end{gather*}
for all $\boldsymbol{\eta } \in L^2(0,T;\boldsymbol{V})$. Therefore, 
the estimates \eqref{m2} for $\mu _{\lambda,\varepsilon}$ and 
$\mu _{\Gamma, \lambda, \varepsilon}$ imply \eqref{m2e}. 
Thus, using \eqref{m2d}--\eqref{m2e}, we show \eqref{m2c}. 
\hfill $\Box$ 

\bigskip
We have obtained the uniform estimates \eqref{m1}, \eqref{m2}, \eqref{m2b}, \eqref{m2c} provided in 
Lemmas~3.1 and 3.2 independent of $\varepsilon $ and $\lambda \in (0,1]$, actually 
$C^* \searrow 1$ as $\varepsilon \searrow 0$. 
They can be used throughout this paper by considering the 
limiting procedure $\varepsilon \searrow 0$ for each fixed $\lambda \in (0,1]$, 
and next the 
limiting procedure $\lambda \searrow  0$.

\section{Proof of the main theorem}
\setcounter{equation}{0}

In this section, we prove Theorem 2.1.

\subsection{Additional uniform estimate independent of $\varepsilon \in (0,1]$}

In this subsection, we obtain additional uniform estimates independent of $\varepsilon \in (0,1]$, 
these estimates may depend on $\lambda \in (0,1]$. 
Therefore, we use them only in the nest subsection to consider 
the limiting procedure as $\varepsilon \searrow 0$.

\paragraph{Lemma 4.1.}
{\it There exist two positive constants $M_3(\lambda )$ and $M_4(\lambda )$, 
depending on $\lambda \in (0,1]$ but 
independent of $\varepsilon \in (0,1]$, such that}
\begin{gather}
	\bigl| \beta _\lambda (u_{\lambda,\varepsilon}) \bigr|_{L^\infty (0,T;H)} 
	\le M_3(\lambda ), \quad 
	\bigl| \beta _{\Gamma, \lambda} (u_{\Gamma, \lambda,\varepsilon}) \bigr|_{L^\infty (0,T;H_\Gamma )} 
	\le M_3(\lambda ),
	\label{m3}
	\\
	|\mu _{\Gamma, \lambda, \varepsilon}|_{L^2(0,T;V_\Gamma)}
	\le M_4(\lambda ), \quad 
	|\mu _{\lambda, \varepsilon}|_{L^2(0,T;V)} 
	\le M_4(\lambda ).
	\label{m4} 
\end{gather}

\paragraph{Proof.}
From the {L}ipschitz continuity of $\beta _\lambda$ and $\beta _{\Gamma, \lambda}$, we see from 
\eqref{m1} that there exists a positive constant $M_3(\lambda )$, 
which is proportional to the {L}ipschitz constants $1/\lambda $ of $\beta _\lambda $ and 
$1/(\lambda \varrho )$ of $\beta _{\Gamma, \lambda}$, 
such that \eqref{m3} holds. 
Next, let us point out the variational equality, deduced from \eqref{apDef2} and 
\eqref{apDef4}:
\begin{align}
	& \varepsilon \bigl( \mu _{\lambda, \varepsilon}(s),z \bigr)_{H} 
	+ \bigl( \mu _{\Gamma, \lambda, \varepsilon}(s),z_\Gamma \bigr)_{H_\Gamma }
	\nonumber \\
	& \quad = \tau \bigl( \partial _s u_{\lambda, \varepsilon}(s),z \bigr)_{H} 
	+ {\revis \bigl( \nabla u_{\lambda, \varepsilon}(s), \nabla z \bigr)_{H} }
	+ \bigl( \beta_\lambda \bigl( u_{\lambda, \varepsilon}(s) \bigr)
	+ \pi \bigl( u_{\lambda, \varepsilon}(s) \bigr)-f_\varepsilon (s),z \bigr)_{H} 
	\nonumber \\
	& \quad \quad {}
	+ \varepsilon \bigl( \partial _s u_{\Gamma, \lambda, \varepsilon} (s), z_\Gamma \bigr)_{H_\Gamma }
	+ {\revis \bigl( \nabla_\Gamma  u_{\Gamma ,\lambda, \varepsilon}(s), \nabla_\Gamma  z_\Gamma
	 \bigr)_{H_\Gamma } }
	\nonumber \\
	& \quad \quad {}
	+ \bigl( \beta_{\Gamma, \lambda }\bigl( u_{\Gamma, \lambda, \varepsilon}(s) \bigr)
	+ \pi_\Gamma  \bigl( u_{\Gamma , \lambda, \varepsilon}(s) \bigr)-f_{\Gamma, \varepsilon}(s),z_\Gamma 
	\bigr)_{H_\Gamma } 
	\label{vari2}
\end{align}
for all $\boldsymbol{z} \in \boldsymbol{V}$, for a.a.\ $s \in (0,T)$. 
Taking now $\boldsymbol{z}:=\boldsymbol{\mu }_{\lambda, \varepsilon }(s)$ 
in \eqref{vari2}, integrating with respect to time, 
we infer with the help of {P}oincar\'e--Wirtinger inequality \eqref{poin1}
\begin{align*} 
	& \varepsilon \int_{0}^{t} 
	\bigl|
	\mu_{\lambda, \varepsilon}(s) 
	\bigr| _H^2 ds 
	+ \int_{0}^{t}
	\bigl|
	\mu_{\Gamma, \lambda, \varepsilon}(s) 
	\bigr| _{H_\Gamma }^2 ds 
	\nonumber \\
	& \le 
	C_P \tau 
	|
	\partial _s u_{\lambda, \varepsilon}
	| _{L^2(0,T;H)}
	\left\{ 
	{\revis 
	|
	\nabla \mu_{\lambda, \varepsilon}
	| _{L^2(0,T;H)}^2
	}
	+
	|
	\mu_{\Gamma, \lambda, \varepsilon}
	| _{L^2(0,T;H_\Gamma) }^2
	\right\}^{1/2} 
	\nonumber \\
	& \quad {}
	+ 
	{\revis 
	|
	\nabla u_{\lambda, \varepsilon}
	| _{L^2(0,T;H)} 
	|
	\nabla \mu _{\lambda, \varepsilon}
	| _{L^2(0,T;H)}
	}
	\nonumber \\
	& \quad {} 
	+ C_P 
	\bigl|
	\beta_\lambda ( u_{\lambda, \varepsilon})
	+ \pi ( u_{\lambda, \varepsilon})-f_\varepsilon 
	\bigr| _{L^2(0,T;H)}
	\left\{ 
	{\revis 
	|
	\nabla \mu _{\lambda, \varepsilon} 
	| _{L^2(0,T;H)}^2 
	}
	+ 
	|
	\mu _{\Gamma, \lambda, \varepsilon} 
	| _{L^2(0,T;H_{\Gamma })}^2 
	\right\}^{1/2} 
	\nonumber \\
	& \quad {} 
	+
	\varepsilon 
	|
	\partial _s u_{\Gamma, \lambda, \varepsilon}
	| _{L^2(0,T;H_\Gamma )} 
	|
	\mu_{\Gamma, \lambda, \varepsilon}
	| _{L^2(0,T;H_\Gamma )}
	+ 
	{\revis 
	|
	\nabla _\Gamma u_{\Gamma, \lambda, \varepsilon}
	| _{L^2(0,T;H_\Gamma )}
	|
	\nabla _\Gamma \mu _{\Gamma, \lambda, \varepsilon} 
	| _{L^2(0,T;H_{\Gamma })}
	}
	\nonumber \\
	& \quad {} 
	+ 
	\bigl|
	\beta_{\Gamma, \lambda} ( u_{\Gamma, \lambda, \varepsilon})
	+ \pi_\Gamma  ( u_{\Gamma, \lambda, \varepsilon})-f_{\Gamma, \varepsilon}
	\bigr| _{L^2(0,T;H_\Gamma )} 
	|
	\mu _{\Gamma, \lambda, \varepsilon} 
	| _{L^2(0,T;H_\Gamma )}
\end{align*}
for all $t \in [0,T]$. 
Then, using the {Y}oung inequality 
along with \eqref{apf1}, \eqref{m1}, \eqref{m2} and \eqref{m3}, 
we deduce the uniform estimate of 
$\{ \mu_{\Gamma, \lambda, \varepsilon} \}_{\varepsilon \in (0,1]}$ in 
$L^2(0,T;V_\Gamma )$. 
Next, 
with the help of
the {P}oincar\'e--{W}irtinger inequality \eqref{poin1} again, 
we deduce the uniform estimate of 
$\{ \mu_{\lambda, \varepsilon} \}_{\varepsilon \in (0,1]}$ in 
$L^2(0,T;H)$. Thus, combining the resultant with \eqref{m2}, 
we see that there exists a positive constant $M_4(\lambda )$, depending on $M_3(\lambda) $ 
independent of $\varepsilon \in (0,1]$, such that \eqref{m4} holds. 
\hfill $\Box$

\paragraph{Lemma 4.2.}
{\it There exists a positive constant $M_5(\lambda )$, 
depending on $\lambda \in (0,1]$ but
independent of $\varepsilon \in (0,1]$, such that}
\begin{equation}
	|u_{\lambda, \varepsilon} |_{L^2(0,T;W)} 
	+
	|u_{\Gamma, \lambda, \varepsilon} |_{L^2(0,T;W_\Gamma)} 
	\le M_5(\lambda ).
	\label{m5}
\end{equation}

\paragraph{Proof.} 
We can compare the terms in \eqref{apDef2} and 
conclude that $\{ |\Delta u_{\lambda, \varepsilon} |_{L^2(0,T;H)} \}_{\varepsilon  \in (0,1]}$
is uniformly bounded.  
Hence, applying the theory of the elliptic regularity 
(see, e.g., \cite[Theorem~3.2, p.~1.79]{BG87}), we have that 
\begin{gather*}
	|u_{\lambda, \varepsilon} |_{L^2(0,T;H^{3/2}(\Omega ))}
	\le \tilde{M}_5(\lambda )\end{gather*}
for some positive constant $\tilde{M}_5(\lambda )$, 
and owing to {\takeshi both the uniform bounds}, we see that
\begin{gather*}
	|\partial_{\boldsymbol{\nu }} u_{\lambda, \varepsilon} |_{L^2(0,T;H_\Gamma )} 
	\le \tilde{M}_5(\lambda ).
\end{gather*}
Next, by comparison in \eqref{apDef4}, 
$\{ |\Delta_\Gamma  u_{\Gamma, \lambda, \varepsilon} |_{L^2(0,T;H_\Gamma)} \}_{\varepsilon  \in (0,1]}$
is uniformly bounded and consequently (see, e.g., \cite[Section~4.2]{Gri09}),
\begin{gather*} 
	| u_{\Gamma, \lambda, \varepsilon}|_{L^2(0,T;W_\Gamma )} 
	 \le \left( | u_{\Gamma, \lambda, \varepsilon}|_{L^2(0,T;V_\Gamma)}^2+
	|\Delta _\Gamma  u_{\Gamma, \lambda, \varepsilon}|_{L^2(0,T;H_\Gamma)}^2 \right)^{1/2} 
	 \le M_5(\lambda ).
\end{gather*} 
for some constant $M_5(\lambda )$. Then, using the theory of the 
elliptic regularity (see, e.g., \cite[Theorem~3.2, p.~1.79]{BG87}), 
we get the conclusion 
\eqref{m5}. 
\hfill $\Box$

\subsection{Passage to the limit as $\varepsilon \searrow  0$}

In this subsection, 
we pass to the limit in the approximating problem as $\varepsilon \searrow  0$. Indeed, 
owing to the estimates stated in Lemmas~3.1, 3.2, 4.1 and 4.2, there 
exist a subsequence of $\varepsilon $ (not relabeled) and some 
limit functions
$u_\lambda, u_{\Gamma, \lambda}$, $\mu _\lambda $ and $\mu _{\Gamma, \lambda}$
 such that 
\begin{gather} 
	u_{\lambda, \varepsilon } \to u_\lambda 
	\quad \mbox{weakly star in } 
	H^1(0,T;H) 
	\cap 
	L^\infty (0,T;V) 
	\cap 
	L^2 (0,T;W), 
	\label{lwcu}
	\\
	\mu _{\lambda, \varepsilon} 
	\to \mu _\lambda 
	\quad \mbox{weakly in } 
	L^2(0,T;V),
	\label{lwcmu}
	\\
	u_{\Gamma, \lambda, \varepsilon } \to u_{\Gamma, \lambda }
	\quad \mbox{weakly star in } 
	H^1(0,T;V_\Gamma ^*) 
	\cap 
	L^\infty (0,T;V_\Gamma ) 
	\cap 
	L^2 (0,T;W_\Gamma ), 
	\label{lwcug}
	\\
	\varepsilon \partial _t u_{\Gamma, \lambda, \varepsilon } \to 0
	\quad \mbox{strongly in } 
	L^2 (0,T;H_\Gamma ), 
	\label{lscug0}
	\\
	\mu _{\Gamma, \lambda, \varepsilon} 
	\to \mu _{\Gamma, \lambda }
	\quad \mbox{weakly in } 
	L^2(0,T;V_\Gamma )
	\label{lwcmug}
\end{gather} 
as $\varepsilon \searrow  0$. 
From \eqref{lwcu} and \eqref{lwcug}, using well-known compactness results 
(see, e.g., \cite[Section~8, Corollary~4]{Sim87}) 
we obtain
\begin{gather} 
	u_{\lambda, \varepsilon} \to u_\lambda  \quad \mbox{strongly in } 
	C\bigl( [0,T];H \bigr) \cap L^2 (0,T;V),
	\label{lscu}
	\\
	u_{\Gamma, \lambda, \varepsilon} \to u_{\Gamma, \lambda}  \quad \mbox{strongly in } 
	C\bigl( [0,T];H_\Gamma  \bigr) \cap L^2 (0,T;V_\Gamma )
	\label{lscug}
\end{gather}
{\takeshi as $\varepsilon \searrow  0$. Then, \eqref{lwcu}, \eqref{lscu} and \eqref{lscug} imply that }
\begin{gather} 
	\varepsilon \partial _t u_{\lambda, \varepsilon } \to 0
	\quad \mbox{strongly in } 
	L^2 (0,T;H), 
	\label{lscug1}
	\\
	\beta_\lambda (u_{\lambda, \varepsilon}) 
	\to \beta_\lambda (u_{\lambda}) \quad \mbox{strongly in } 
	C\bigl( [0,T];H \bigr),
	\label{lscb}
	\\
	\beta_{\Gamma, \lambda} (u_{\Gamma, \lambda, \varepsilon}) 
	\to \beta_{\Gamma, \lambda} (u_{\Gamma, \lambda}) 
	\quad \mbox{strongly in } 
	C\bigl( [0,T];H_\Gamma \bigr)
	\label{lscbg}
	\\
	\pi(u_{\lambda, \varepsilon}) \to \pi (u_{\lambda}) \quad \mbox{strongly in } 
	C\bigl( [0,T];H \bigr),
	\label{lscp}
	\\
	\pi_{\Gamma} (u_{\Gamma, \lambda, \varepsilon}) \to \pi_{\Gamma} (u_{\Gamma, \lambda}) 
	\quad \mbox{strongly in } 
	C\bigl( [0,T];H_\Gamma \bigr)
	\label{lscpg}
\end{gather}
as $\varepsilon \searrow  0$. 
We point out that \eqref{apu}, \eqref{lscu} and \eqref{lscug} 
entail that
\begin{equation}
	u_\lambda (0)=u_0 \quad \mbox{a.e.\ in } \Omega, 
	\quad 
	u_{\Gamma, \lambda}(0) = u_{0\Gamma } 
	\quad \mbox{a.e.\ on } \Gamma.
	\label{lini}
\end{equation}

From \eqref{apDef1} and 
\eqref{apDef3} it follows that
\begin{align*}
	& \varepsilon \bigl( \partial _s u _{\lambda, \varepsilon}(s),z \bigr)_{H} 
	+ \bigl( \partial _s u_{\Gamma, \lambda, \varepsilon}(s),z_\Gamma \bigr)_{H_\Gamma }
	+ 
	{\revis 
	\bigl( \nabla \mu _{\lambda, \varepsilon}(s), \nabla z \bigr)_{H} 
	}
	+ 
	{\revis 
	\bigl( \nabla_\Gamma \mu _{\Gamma ,\lambda, \varepsilon}(s), \nabla_\Gamma  z_\Gamma
	 \bigr)_{H_\Gamma } 
	 }
	= 0
\end{align*}
for all $\boldsymbol{z} \in \boldsymbol{V}$, for a.a.\ $s \in (0,T)$. 
Then using 
\eqref{lwcu}--\eqref{lscpg}
we can pass to the limit in this variational equality and in \eqref{vari2} obtaining
\begin{gather}
	\bigl( \partial _s u_{\Gamma, \lambda}(s),z_\Gamma \bigr)_{H_\Gamma }
	+ 
	{\revis 
	\bigl( \nabla \mu _{\lambda}(s), \nabla z \bigr)_{H} 
	}
	+ 
	{\revis 
	\bigl( \nabla_\Gamma \mu _{\Gamma ,\lambda}(s), \nabla_\Gamma  z_\Gamma
	 \bigr)_{H_\Gamma } 
	 }
	= 0,
	\label{vari3}
	\\
	\bigl( \mu _{\Gamma, \lambda}(s),z_\Gamma \bigr)_{H_\Gamma }
	= \tau \bigl( \partial _t u_{\lambda}(s),z \bigr)_{H} 
	+ 
	{\revis 
	\bigl( \nabla u_{\lambda}(s), \nabla z \bigr)_{H} 
	}
	+ \bigl( \beta_\lambda \bigl( u_{\lambda}(s) \bigr)
	+ \pi \bigl( u_{\lambda}(s) \bigr)-f(s),z \bigr)_{H} 
	\nonumber \\
	{} 
	+ 
	{\revis 
	\bigl( \nabla_\Gamma  u_{\Gamma ,\lambda}(s), \nabla_\Gamma  z_\Gamma
	 \bigr)_{H_\Gamma } 
	 }
	+ \bigl( \beta_{\Gamma, \lambda }\bigl( u_{\Gamma, \lambda}(s) \bigr)
	+ \pi_\Gamma  \bigl( u_{\Gamma , \lambda}(s) \bigr)-f_\Gamma (s),z_\Gamma 
	\bigr)_{H_\Gamma } 
	\label{vari4}
\end{gather}
for all $\boldsymbol{z} \in \boldsymbol{V}$. Moreover, from a priori estimates obtained in 
Lemmas~3.1 and 3.2, we have 
\begin{gather} 
	|u_\lambda|_{L^\infty (0,T;V)} \le M_1, \quad 
	|u_{\Gamma, \lambda}|_{L^\infty (0,T;H_\Gamma )} \le c_1^{1/2} M_1,
	\label{m1n}\\
	|\partial _t u_\lambda|_{L^2(0,T;H)}
	+
	{\revis 
	|\nabla _\Gamma u_{\Gamma, \lambda}|_{L^\infty (0,T;H_\Gamma)}
	}
	+
	\bigl|
	\widehat{\beta }_{\lambda } ( u_\lambda) 
	\bigr|_{L^\infty (0,T;L^1(\Omega ))}
	\nonumber \\
	\quad {}+ 
	\bigl|
	\widehat{\beta }_{\Gamma, \lambda } ( u_{\Gamma, \lambda} ) 
	\bigr|_{L^\infty (0,T;L^1(\Gamma))}
	+
	{\revis 
	|\nabla \mu _{\lambda}|_{L^2(0,T;H)}
	}
	+
	{\revis 
	|\nabla_\Gamma  \mu _{\Gamma, \lambda}|_{L^2(0,T;H_{\Gamma })}
	}
	\le M_2,
	\label{m2n}\\
	| \partial_t u_{\Gamma, \lambda}
	|_{L^2(0,T;V_\Gamma ^*)} 
	\le 2M_2,
	\label{m2cn}
\end{gather}
because $C^* \searrow 1$ as $\varepsilon \searrow 0$. 
As a remark, taking $z=1$, $z_\Gamma =1$ 
in \eqref{vari3} and using \eqref{lscug}, \eqref{lini}, we obtain that 
\begin{equation}
	\int_{\Gamma }^{} u_{\lambda }(t) d\Gamma = \int_{\Gamma }^{} u_{0\Gamma }d\Gamma 
	\quad \mbox{for all~} t \in [0,T].
	\label{vpl}
\end{equation}

\subsection{Proof of Theorem 2.1.}

In this subsection, we prove the main theorem. 
To do so, we are going to produce estimates 
independent of $\lambda $ and passing to the limit as 
$\lambda \searrow  0$. 
The point of emphasis is the effective usage of the mean value zero function.

\paragraph{Lemma 4.3} 
{\it There exists a positive constant $M_6$, independent of $\lambda \in (0,1]$ such that} 
\begin{equation}
	\bigl| 
	\beta_\lambda ( u_{\lambda})
	\bigr|_{L^2(0,T;L^1(\Omega))}
	+ 
	\bigl| 
	\beta_{\Gamma,  \lambda} ( u_{\Gamma, \lambda})
	\bigr|_{L^2(0,T;L^1(\Gamma))}
	\le M_6.
	\label{m6}
\end{equation}

\paragraph{Proof.} 
Recall \eqref{mgamma} and take 
\begin{equation*}
	z:=u_{\lambda }(s)-m_\Gamma, \quad 
	z_\Gamma :=u_{\Gamma,\lambda}(s)-m_\Gamma 
\end{equation*}
in \eqref{vari4}. Then we have 
\begin{gather}
	\tau \bigl( \partial _t u_{\lambda}(s),u_{\lambda }(s)-m_\Gamma \bigr)_{H} 
	+
	{\revis 
	 \bigl| \nabla u_{\lambda}(s) \bigr|_{H}^2
	 }
	+ \bigl( \beta_\lambda \bigl( u_{\lambda}(s) \bigr), u_{\lambda }(s)-m_\Gamma \bigr)_{H} 
	\nonumber \\
	{}+ \bigl( \pi \bigl( u_{\lambda}(s) \bigr)-f(s), u_{\lambda }(s)-m_\Gamma \bigr)_{H} 
	+ 
	{\revis 
	\bigl| \nabla_\Gamma  u_{\Gamma ,\lambda}(s) \bigr|_{H_\Gamma }^2
	}
	+ \bigl( \beta_{\Gamma, \lambda }\bigl( u_{\Gamma, \lambda}(s) \bigr), 
	u_{\Gamma,\lambda}(s)-m_\Gamma 
	\bigr)_{H_\Gamma } 
	\nonumber \\
	{}+ \bigl( \pi_\Gamma \bigl( u_{\Gamma , \lambda}(s) \bigr)-f_\Gamma (s),u_{\Gamma,\lambda}(s)-m_\Gamma 
	\bigr)_{H_\Gamma } 
	= \bigl( \mu _{\Gamma, \lambda}(s),u_{\Gamma,\lambda}(s)-m_\Gamma  \bigr)_{H_\Gamma }.
	\label{vari5}
\end{gather}
Let now $(y,y_\Gamma )$ be the solution of the following problem
\begin{gather}
	\int_{\Omega }^{} \nabla y\cdot \nabla z dx 
	+ \int_{\Gamma }^{} \nabla _\Gamma y_\Gamma \cdot \nabla _\Gamma z_\Gamma d\Gamma 
	= \int_{\Gamma }^{} (u_{\Gamma,\lambda}-m_\Gamma ) z_\Gamma d\Gamma 
	\quad \mbox{for all } \boldsymbol{z} \in \boldsymbol{V}, 
	\label{aux}
	\\
	y_\Gamma =y_{|_\Gamma } \quad \mbox{a.e.\ on } \Gamma, 
	\quad 
	\int_{\Gamma }^{} y_\Gamma d\Gamma =0 \quad \mbox{a.e.\ in } (0,T). 
	\nonumber 
\end{gather}
This problem has one and only one solution, see \cite[Appendix, Lemma~A]{CF15b} and 
repeat the proof with the different condition. Moreover, due to \eqref{poin1} and \eqref{poin2}, 
there exists a positive constant $\tilde{M}_6$, independent of $\lambda \in (0,1]$ such that
\begin{equation} 
	\bigl| y(s) \bigr|_{V}^2 + \bigl| y_\Gamma (s) \bigr|_{V_\Gamma }^2 
	\le \tilde{M}_6 \bigl| u_{\Gamma, \lambda  }(s)-m_\Gamma \bigr|_{H_\Gamma }^2 
	\label{m6t}
\end{equation}
for a.a.\ $s \in (0,T)$. 
Taking $z:=\mu _\lambda (s)$, $z_\Gamma :=\mu _{\Gamma, \lambda }$ in \eqref{aux}, 
we see that the right hand side of \eqref{aux} is equal to 
\begin{equation*}
	{\revis
	\bigl( \nabla y(s), \nabla \mu_\lambda (s) \bigr)_{H}
	}
	+ 
	{\revis 
	\bigl( \nabla _\Gamma y_\Gamma(s), \nabla _\Gamma \mu _{\Gamma, \lambda }(s) \bigr)_{H_\Gamma }
	} 
\end{equation*} 
and consequently we use \eqref{vari3} with $z:=y(s)$, $z_\Gamma :=y_\Gamma (s)$ to conclude that 
\begin{equation*}
	\bigl( \mu _{\Gamma ,\lambda}(s), u_{\Gamma, \lambda }(s) - m_\Gamma 
	 \bigr)_{H_\Gamma } 
	= - \bigl( \partial _s u_{\Gamma, \lambda}(s),y_\Gamma(s) \bigr)_{H_\Gamma }. 
\end{equation*}
Then, from \eqref{vari5} and the properties that 
\begin{gather*}
	\beta _\lambda  (r)(r-m_\Gamma ) \ge \delta _0 
	\bigl | \beta _\lambda(r) \bigr | -c_3,
	\quad 
	\beta _{\Gamma, \lambda} (r)(r-m_\Gamma ) \ge \delta _0 
	\bigl | \beta _{\Gamma, \lambda} (r) \bigr | -c_3
\end{gather*}
for all $r \in \mathbb{R}$, $\lambda \in (0,1]$ and some positive constants 
$\delta _0$ and $c_3$ which are provided from \cite[Section~5]{GMS09} with 
the assumptions 
$D(\beta _\Gamma ) \subseteq D(\beta )$ of {\rm (A5)} and 
$m_\Gamma  \in {\rm int} D(\beta _\Gamma )$ of {\rm (A6)}, we deduce
\begin{align*}
	& 
	{\revis 
	\bigl| \nabla u_{\lambda}(s) \bigr|_{H}^2 
	}
	+ 
	\delta _0 \bigl| \beta_\lambda \bigl( u_{\lambda}(s) \bigr)
	\bigr|_{L^1(\Omega)}
	+ 
	{\revis 
	\bigl| \nabla_\Gamma  u_{\Gamma ,\lambda}(s) \bigr|_{H_\Gamma }^2 
	}
	+
	\delta _0  \bigl| \beta_{\Gamma, \lambda} \bigl( u_{\Gamma, \lambda}(s) \bigr)
	\bigr|_{L^1(\Gamma )} 
	\nonumber \\
	& 
	\le c_3T\bigl( |\Omega |+|\Gamma | \bigr)+ 
	\bigl| \tau  \partial _t u_{\lambda}(s) + \pi \bigl( u_\lambda(s) \bigr)-f(s)
	\bigr|_{H} 
	\bigl| 
	u_{\lambda }(s)-m_\Gamma 
	\bigr|_{H}
	\nonumber \\
	& \quad {} + \bigl| \pi_\Gamma \bigl( u_{\Gamma , \lambda}(s) \bigr)-f_\Gamma(s)
	\bigr|_{H_\Gamma}
	\bigl| u_{\Gamma,\lambda}(s)-m_\Gamma 
	\bigr|_{H_\Gamma} 
	+ 
	\bigl| \partial _s u_{\Gamma, \lambda}(s) \bigr|_{V_\Gamma^*} 
	\bigl| y_\Gamma(s) \bigr|_{V_\Gamma }.
\end{align*}
Now, squaring both sides we obtain 
\begin{align}
	& \left( 
	\delta _0 \bigl| \beta_\lambda \bigl( u_{\lambda}(s) \bigr)
	\bigr|_{L^1(\Omega)}
	+ 
	\delta _0  \bigl| \beta_{\Gamma, \lambda} \bigl( u_{\Gamma, \lambda}(s) \bigr)
	\bigr|_{L^1(\Gamma )} 
	\right)^2
	\nonumber \\
	& 
	\le 4c_3^2T^2\bigl( |\Omega |+|\Gamma | \bigr)^2+ 
	12 \left( 
	\tau ^2 \bigl| \partial _t u_{\lambda}(s)\bigr|^2_H + 
	\bigl| \pi \bigl( u_\lambda(s) \bigr) \bigr|_H ^2 + 
	\bigl| f(s)
	\bigr|_{H}^2
	\right)  
	\bigl| 
	u_{\lambda }(s)-m_\Gamma 
	\bigr|_{H}^2
	\nonumber \\
	& \quad {} + 8 \left( \bigl| \pi_\Gamma \bigl( u_{\Gamma , \lambda}(s) \bigr) \bigl|_{H_\Gamma }^2 
	+ \bigl| f_\Gamma(s) \bigr|_{H_\Gamma}^2 
	\right) 
	\bigl| u_{\Gamma,\lambda}(s)-m_\Gamma 
	\bigr|_{H_\Gamma}^2
	+ 
	\bigl| \partial _s u_{\Gamma, \lambda}(s) \bigr|_{V_\Gamma^*}^2 
	\bigl| y_\Gamma(s) \bigr|_{V_\Gamma }^2
	\label{l2l1}
\end{align}
for a.a.\ $s \in (0,T)$
Here, by virtue of \eqref{m1n}, \eqref{poin2}, \eqref{m2n} and \eqref{m6t}
\begin{gather*}
	| 
	u_{\lambda }-m_\Gamma 
	|_{L^\infty (0,T;H)}^2 \le 2M_1^2+2m_\Gamma |\Omega |, 
	\\
	| 
	u_{\Gamma,\lambda}-m_\Gamma 
	|_{L^\infty (0,T;H_\Gamma)}^2
	\le C_P M_2^2,
	\\
	| 
	y_\Gamma 
	|_{L^\infty (0,T;V_\Gamma) }^2 \le C_P M_2^2 \tilde{M}_6
\end{gather*}
because $\int_{\Gamma }^{}(u_{\Gamma,\lambda}(s)-m_\Gamma )d\Gamma =0$ for a.a.\ $s \in (0,T)$. 
Therefore, integrate \eqref{l2l1} over $[0,T]$ with respect to time.  
Then 
the right hand side can be bounded due to \eqref{m1n}, \eqref{m2n} and \eqref{m2cn}. 
Thus, there exists a positive constant $M_6$, independent of $\lambda \in (0,1]$, 
such that 
\eqref{m6} holds. \hfill $\Box$ 

\bigskip
Put 
\begin{equation*}
	\omega_\lambda (t) 
	:= \frac{1}{|\Gamma |}\int_{\Gamma }^{} 
	\mu _{\Gamma, \lambda }(t) dt 
\end{equation*}
for a.a.\ $t \in (0,T)$. Then we obtain the following estimate.

\paragraph{Lemma 4.4.} 
{\it There exist two positive constants $M_7$ and $M_8$ 
independent of $\lambda \in (0,1]$ such that} 
\begin{gather} 
	|
	\omega _\lambda 
	|_{L^2(0,T)} 
	\le M_7,
	\label{m7}
	\\
	| 
	\mu _{\Gamma, \lambda} 
	|_{L^2(0,T;V_\Gamma )}
	\le M_8, \quad 
	| 
	\mu _\lambda 
	|_{L^2(0,T;V)}
	\le M_8.
	\label{m8}
\end{gather}

\paragraph{Proof.}
Taking $z:=1/|\Gamma |$ and $z_\Gamma :=1/|\Gamma |$ in 
\eqref{vari4} we obtain 
\begin{align*}
	\bigl| 
	\omega _\lambda (t) 
	\bigr|^2
	& = \frac{7}{|\Gamma |^2} \left\{ 
	\tau^2 
	\bigl| \partial _t u_{\lambda }(t) \bigr|_{L^1(\Omega)}^2
	+ \bigl| 
	\beta _\lambda 
	\bigl( 
	u_{\lambda }(t)
	\bigr) 
	\bigr|_{L^1(\Omega )}^2
	+
	\bigl|
	\pi \bigl( u_{\lambda }(t)
	\bigr) 
	\bigr|_{L^1(\Omega)}^2
	+ 
	\bigl|
	f_\lambda (t) 
	\bigr|_{L^1(\Omega)}^2
	\right.
	\nonumber \\
	& \quad 
	\left. 
	{}
	+ 
	\bigl| 
	\beta _{\Gamma, \lambda} 
	\bigl( 
	u_{\Gamma, \lambda }(t)
	\bigr) 
	\bigr|_{L^1(\Gamma)}^2
	+
	\bigl| 
	\pi_\Gamma  \bigl( u_{\Gamma, \lambda }(t)
	\bigr) 
	\bigr|_{L^1(\Gamma)}^2
	+ \bigl| 
	f_{\Gamma, \lambda} (t) 
	\bigr|_{L^1(\Gamma)}^2
	\right\}
\end{align*}
for a.a.\ $t \in (0,T)$, that is, 
there exists a positive constant $M_7$
independent of $\lambda \in (0,1]$ such that
\eqref{m7} holds. 
Next, using \eqref{poin2} we obtain 
\begin{align}
	\bigl|
	\mu _{\Gamma, \lambda} (t)
	\bigr|_{H_\Gamma }^2
	& \le 2
	\bigl|
	\mu _{\Gamma, \lambda} (t)-\omega _\lambda (t)
	\bigr|_{H_\Gamma }^2
	+2 \bigl| \omega _\lambda (t) \bigr|^2_{H_\Gamma }
	\nonumber \\
	& \le {\revis 
	2C_P \bigl|
	\nabla _\Gamma \mu _{\Gamma, \lambda} (t)
	\bigr|_{H_\Gamma}^2
	}
	+2|\Gamma | \bigl| \omega _\lambda (t) \bigr|^2
\end{align}
for a.a.\ $t\in (0,T)$. 
Thus, \eqref{m2n} and \eqref{m7} imply the 
first estimate of \eqref{m8} with some positive constant $M_8$ independent of $\lambda \in (0,1]$. 
Moreover, by using \eqref{poin1} with the above, \eqref{m2n} ensures the 
validity of the second estimate in \eqref{m8}. 
\hfill $\Box$. 

\bigskip
\paragraph{Lemma 4.5.} 
{\it There exist a positive constant $M_9$
independent of $\lambda \in (0,1]$ such that} 
\begin{gather} 
	\bigl|
	\beta _\lambda 
	(
	u_\lambda
	)
	\bigr|_{L^2(0,T;H)} 
	+
	\bigl|
	\beta _{\lambda }
	(
	u_{\Gamma, \lambda}
	)
	\bigr|_{L^2(0,T;H_\Gamma )} 
	\le M_9,
	\label{m9a}
	\\
	| 
	\Delta u_{\lambda }
	|_{L^2(0,T;H)}
	+ 
	| 
	u_{\lambda} 
	|_{L^2(0,T;H^{3/2}(\Omega) )}
	+ 
	| 
	\partial _{\boldsymbol{\nu }}u_{\lambda} 
	|_{L^2(0,T;H_\Gamma )}
	\le M_9,
	\label{m9b} \\
	\bigl|
	\beta _{\Gamma, \lambda }
	(
	u_{\Gamma, \lambda}
	)
	\bigr|_{L^2(0,T;H_\Gamma )} 
	\le M_9,
	\label{m9c} \\
	|
	u_\lambda
	|_{L^2(0,T;W)} 
	+
	|
	u_{\Gamma, \lambda}
	|_{L^2(0,T;W_\Gamma )} 
	\le M_9.
	\label{m9d}
\end{gather}

The estimates \eqref{m9a}--\eqref{m9b} take advantage from writing 
\eqref{vari4} as the combination of the equations
\begin{gather}
	\tau \partial _t u_\lambda 
	-
	\Delta u_\lambda 
	+ 
	\beta _\lambda (u_\lambda )
	+
	\pi (u_\lambda )
	=f
	\quad \mbox{a.e.\ in } Q,
	\label{lequ1}
	\\
	\mu _{\Gamma,\lambda}
	=\partial _{\boldsymbol{\nu }}
	u_\lambda 
	-\Delta _\Gamma 
	u_{\Gamma, \lambda}
	+\beta _{\Gamma, \lambda} (u_{\Gamma, \lambda  })
	+\pi _\Gamma (u_{\Gamma, \lambda  })
	- f _{\Gamma}
	\quad \mbox{a.e.\ on } \Sigma,
	\label{lequ2}
\end{gather}
which are rigorous due to the regularity of $u_\lambda$ and 
$u_{\Gamma, \lambda }$ stated in \eqref{m5}. 
Recalling \cite[Lemmas~4.4, 4.5]{CC13}, we observe that {\revis the 
proof} is essentially the same as in these lemmas. 
Therefore we omit the details for \eqref{m9a}--\eqref{m9d}. 

\medskip
We have collected all information which enables us to pass to the limit 
as $\lambda \searrow  0$. 

\medskip
\paragraph{Proof of Theorem 2.1.} 
Thanks to \eqref{m1n}--\eqref{m2cn}, \eqref{m8}, \eqref{m9a}, \eqref{m9c} and \eqref{m9d}, 
we see that 
there 
exist a subsequence of $\lambda$ (not relabeled) and some 
limit functions
$u, u_{\Gamma}$, $\mu$, $\mu _\Gamma$, $\xi$ and $\xi _\Gamma $
 such that 
\begin{gather} 
	u_\lambda \to u
	\quad \mbox{weakly star in } 
	H^1(0,T;H) 
	\cap 
	L^\infty (0,T;V) 
	\cap 
	L^2 (0,T;W), 
	\label{wcu}
	\\
	\mu _\lambda
	\to \mu
	\quad \mbox{weakly in } 
	L^2(0,T;V),
	\label{wcmu}
	\\
	u_{\Gamma, \lambda} \to u_{\Gamma}
	\quad \mbox{weakly star in } 
	H^1(0,T;V_\Gamma ^*) 
	\cap 
	L^\infty (0,T;V_\Gamma ) 
	\cap 
	L^2 (0,T;W_\Gamma ), 
	\label{wcug}
	\\
	\mu _{\Gamma, \lambda} 
	\to \mu _\Gamma
	\quad \mbox{weakly in } 
	L^2(0,T;V_\Gamma ),
	\label{wcmug}
	\\
	\beta_\lambda (u_\lambda) 
	\to \xi \quad \mbox{weakly in } 
	L^2(0,T;H),
	\label{wcb}
	\\
	\beta_{\Gamma, \lambda} (u_{\Gamma, \lambda}) 
	\to \xi _\Gamma 
	\quad \mbox{weakly in } 
	L^2(0,T;H_\Gamma)
	\label{wcbg}
\end{gather} 
as $\lambda \searrow  0$. 
From \eqref{wcu} and \eqref{wcug}, using well-known compactness results 
(see, e.g., \cite[Section~8, Corollary~4]{Sim87}) 
we obtain
\begin{gather} 
	u_\lambda \to u  \quad \mbox{strongly in } 
	C\bigl( [0,T];H \bigr) \cap L^2 (0,T;V),
	\label{scu}
	\\
	u_{\Gamma, \lambda} \to u_{\Gamma}  \quad \mbox{strongly in } 
	C\bigl( [0,T];H_\Gamma  \bigr) \cap L^2 (0,T;V_\Gamma ),
	\label{scug}
\end{gather}
which imply that 
\begin{gather} 
	\pi(u_\lambda) \to \pi (u) \quad \mbox{strongly in } 
	C\bigl( [0,T];H \bigr),
	\label{scp}
	\\
	\pi_{\Gamma} (u_{\Gamma, \lambda}) \to \pi_{\Gamma} (u_{\Gamma}) 
	\quad \mbox{strongly in } 
	C\bigl( [0,T];H_\Gamma \bigr)
	\label{scpg}
\end{gather}
as $\lambda \searrow  0$ and 
\begin{equation}
	\xi \in \beta (u) 
	\quad \mbox{a.e.\ in } Q, 
	\quad 
	\xi _\Gamma 
	\in \beta_\Gamma (u_\Gamma )
	\quad \mbox{a.e.\ on } \Sigma
	\label{in}
\end{equation}
due to the maximal monotonicity of $\beta $ and $\beta _\Gamma $. 
Finally, we can pass to the limit as {\takeshi $\lambda \searrow 0$} in the variational equality 
\eqref{vari3} as well as in \eqref{lequ1}--\eqref{lequ2} obtaining 
\eqref{Def1}--\eqref{Def3}, and 
in the initial conditions \eqref{lini} obtaining \eqref{Def4}. 
Thus, we arrive at the conclusion. 
\hfill $\Box$

\paragraph{Remark 2.} 
Thanks to the strong convergence \eqref{scug} and the equality \eqref{vpl}, 
we deduce that 
\begin{equation}
	\int_{\Gamma }^{} u(t) d\Gamma = \int_{\Gamma }^{} u_{0\Gamma }d\Gamma 
	\quad \mbox{for all~} t \in [0,T].
	\label{meanae}
\end{equation}

\subsection{Proof of Theorem 2.2.}

\paragraph{Proof of Theorem 2.2.}
Let us note that the solutions 
$\boldsymbol{u}^{(i)}:=(
u^{(i)}, u_\Gamma ^{(i)})$, 
$\boldsymbol{\mu }^{(i)}:=(\mu^{(i)}, \mu _\Gamma ^{(i)})$, 
$\boldsymbol{\xi}^{(i)}:=(\xi^{(i)}, \xi _\Gamma ^{(i)})$ 
of {\rm (P)} satisfy
\begin{gather}
	{\fukao \bigl\langle  \partial _s u_{\Gamma}^{(i)}(s),z_\Gamma \bigr\rangle _{V_\Gamma^*,V_\Gamma } }
	+ {\revis \bigl( \nabla \mu^{(i)}(s), \nabla z \bigr)_{H} }
	+ {\revis \bigl( \nabla_\Gamma \mu _{\Gamma}^{(i)}(s), \nabla_\Gamma  z_\Gamma
	 \bigr)_{H_\Gamma } }
	= 0,
	\label{vari3i}
	\\
	\bigl( \mu _{\Gamma}^{(i)}(s),z_\Gamma \bigr)_{H_\Gamma }
	= \tau \bigl( \partial _t u^{(i)}(s),z \bigr)_{H} 
	+ {\revis \bigl( \nabla u^{(i)}(s), \nabla z \bigr)_{H} }
	+ \bigl( \xi^{(i)}(s) 
	+ \pi \bigl( u^{(i)}(s) \bigr)-f^{(i)}(s),z \bigr)_{H} 
	\nonumber \\
	{} + {\revis \bigl( \nabla_\Gamma  u_{\Gamma}^{(i)}(s), \nabla_\Gamma  z_\Gamma
	 \bigr)_{H_\Gamma } }
	+ \bigl( \xi _{\Gamma}^{(i)}(s) 
	+ \pi_\Gamma  \bigl( u_{\Gamma}^{(i)}(s) \bigr)-f_{\Gamma}^{(i)}(s),z_\Gamma 
	\bigr)_{H_\Gamma } 
	\label{vari4i}
\end{gather}
for all $\boldsymbol{z}:=(z,z_\Gamma ) \in \boldsymbol{V}$ and for $i=1,2$. 
Let us use the notation 
\begin{gather*}
	\bar{u}:=u^{(1)}-u^{(2)}, 
	\quad  
	\bar{u}_\Gamma :=u_\Gamma ^{(1)}-u_\Gamma ^{(2)}, 
	\quad 
	\bar{\mu }:=\mu^{(1)}-\mu^{(2)}, 
	\quad 
	\bar{\mu }_\Gamma :=\mu _\Gamma ^{(1)}-\mu _\Gamma ^{(2)}, 
	\\
	\bar{\xi }:=\xi ^{(1)}-\xi ^{(2)}, 
	\quad 
	\bar{\xi }_\Gamma :=\xi _\Gamma ^{(1)}-\xi _\Gamma ^{(2)}, 
	\quad 
	\bar{f}:=f^{(1)}-f^{(2)}, 
	\quad 
	\bar{f}_\Gamma :=f_\Gamma ^{(1)}-f_\Gamma ^{(2)},
	\\ 
	\bar{u}_0:=u_{0}^{(1)}-u_{0}^{(2)},
	\quad 
	\bar{u}_{0\Gamma}:=u_{0\Gamma }^{(1)}-u_{0\Gamma }^{(2)}. 
\end{gather*}
We take the difference of equations \eqref{vari4i}, test it by 
$(\bar{u}, \bar{u}_\Gamma )$ and integrate over $[0,t]$ with respect to $s$, obtaining 
\begin{align}
	& \frac{\tau }{2} \bigl| \bar{u}(t) \bigr|_H^2 
	+ 
	{\revis 
	\int_{0}^{t} \bigl| \nabla \bar{u}(s) \bigr|_{H}^2 ds 
	}
	+
	{\revis 
	\int_{0}^{t} \bigl| \nabla_\Gamma  \bar{u}_\Gamma (s) \bigr|_{H_\Gamma}^2 ds 
	}
	+
	\int_{0}^{t} \bigl( \bar{\xi }(s), \bar{u}(s) \bigr)_H ds 
	\nonumber \\
	& \quad {}
	+
	\int_{0}^{t} \bigl( \bar{\xi }_\Gamma (s), \bar{u}_\Gamma (s) \bigr)_{H_\Gamma } ds 
	- \int_{0}^{t} \bigl( \bar{\mu }_\Gamma (s), \bar{u}_\Gamma (s) \bigr)_{H_\Gamma } ds 
	\nonumber \\
	& = \frac{\tau }{2} \bigl| \bar{u}_0 \bigr|_H^2 
	-
	\int_{0}^{t} \bigl( \pi \bigl( u^{(1)}(s) \bigr) -\pi \bigl( u^{(2)}(s) \bigr), \bar{u}(s) 
	\bigr)_H ds
	\nonumber \\
	& \quad {}
	- 
	\int_{0}^{t} \bigl( \pi_\Gamma \bigl( u_\Gamma ^{(1)}(s) \bigr) 
	-\pi_\Gamma  \bigl( u_\Gamma ^{(2)}(s) \bigr), \bar{u}_\Gamma (s) 
	\bigr)_{H_\Gamma } ds
	\nonumber \\
	& \quad {}
	+ \int_{0}^{t} \bigl( \bar{f}(s), \bar{u}(s) \bigr)_H ds 
	+ \int_{0}^{t} \bigl( \bar{f}_\Gamma (s), \bar{u}_\Gamma (s) \bigr)_{H_\Gamma } ds 
	\label{last}
\end{align}
for all $t \in [0,T]$. 
We start by discussing the last term in the left hand side of \eqref{last}. Let 
$(\bar{y}, \bar{y}_\Gamma ) \in H^1(0,T;\boldsymbol{V})$ solves the problem 
\begin{equation}
	\int_{\Omega }^{} \nabla \bar{y}\cdot \nabla z dx + 
	\int_{\Gamma }^{} \nabla _\Gamma \bar{y}_\Gamma \cdot \nabla _\Gamma z_\Gamma d\Gamma 
	=\langle \bar{u}_\Gamma , z_\Gamma \rangle _{V^*_\Gamma , V_\Gamma }
	\label{aux1}
\end{equation}
for all $\boldsymbol{z}:=(z,z_\Gamma ) \in \boldsymbol{V}$, with 
\begin{equation}
	\int_{\Gamma }^{} \bar{y}_\Gamma d\Gamma =0,
	\label{aux2}
\end{equation}
a.e.\ in $(0,T)$. 
Now, 
testing \eqref{aux1} by $\boldsymbol{z}=(\bar{\mu }(s), \bar{\mu }_\Gamma(s))$, we 
have 
\begin{equation*}
	- \int_{0}^{t} \bigl( \bar{\mu }_\Gamma (s), \bar{u}_\Gamma (s) \bigr)_{H_\Gamma } ds 
	= 
	-
	{\revis 
	\int_{0}^{t} \bigl( 
	\nabla \bar{y}(s), 
	\nabla \bar{\mu }(s) \bigr)_{H} ds 
	}
	- 
	{\revis 
	\int_{0}^{t} \bigl( 
	\nabla _\Gamma \bar{y}_\Gamma(s),
	\nabla _\Gamma \bar{\mu }_\Gamma(s) \bigr) _{H_\Gamma} ds 
	}
\end{equation*}
for all $t \in [0,T]$. 
On the other hand, considering the difference of 
\eqref{vari3i} and testing then by 
$\boldsymbol{z}=(\bar{y}(s), \bar{y}_\Gamma(s))$, 
we infer that 
\begin{equation*}
	- 
	{\revis 
	\int_{0}^{t} \bigl( \nabla \bar{\mu }(s), \nabla \bar{y}(s) \bigr)_{H} ds 
	}
	- 
	{\revis 
	\int_{0}^{t} \bigl( \nabla _\Gamma \bar{\mu }_\Gamma(s), 
	\nabla _\Gamma \bar{y}_\Gamma(s) \bigr) _{H_\Gamma} ds 
	}
	= \int_{0}^{t} \bigl\langle 
	\partial _s \bar{u}_\Gamma(s) ,\bar{y}_\Gamma(s) \bigr\rangle _{V^*_\Gamma , V_\Gamma }
	ds
\end{equation*}
for all $t \in [0,T]$.
Therefore, we can differentiate \eqref{aux1} with respect to time, then test by $\boldsymbol{z}=(\bar{y}(s),\bar{y}_\Gamma(s))$ and 
integrate the resultant over $[0,t]$ with respect to $s$, obtaining
\begin{align*}
	{\revis 
	\frac{1}{2} \bigl| \nabla \bar{y}(t) \bigr|_{H}^2 
	}
	+ 
	{\revis 
	\frac{1}{2} \bigl| \nabla \bar{y}_\Gamma (t) \bigr|_{H_\Gamma}^2 
	}
	- 
	{\revis 
	\frac{1}{2} \bigl| \nabla \bar{y}(0) \bigr|_{H}^2
	} 
	- 
	{\revis 
	\frac{1}{2} \bigl| \nabla \bar{y}_\Gamma(0) \bigr|_{H_\Gamma}^2 
	}
	= \int_{0}^{t} \bigl\langle 
	\partial _s \bar{u}_\Gamma(s) ,\bar{y}_\Gamma(s) \bigr\rangle _{V^*_\Gamma , V_\Gamma }
	ds
\end{align*}
for all $t\in [0,T]$. 
Thus, we have the contribution 
\begin{equation}
	{\revis 
	\frac{1}{2} \bigl| \nabla \bar{y}(t) \bigr|_{H}^2 
	}
	+ 
	{\revis 
	\frac{1}{2} \bigl| \nabla \bar{y}_\Gamma (t) \bigr|_{H_\Gamma }^2 
	}
	\ge c_4 \bigl| \bar{u}_\Gamma (t) \bigr|_{V_\Gamma ^*}^2
	\label{yt}
\end{equation}
on the left hand side and 
\begin{equation}
	{\revis 
	\frac{1}{2} \bigl| \nabla \bar{y}(0) \bigr|_{H}^2 
	}
	+ 
	{\revis 
	\frac{1}{2} \bigl| \nabla \bar{y}_\Gamma (t) \bigr|_{H_\Gamma}^2 
	}
	\le c_5 | \bar{u}_{0\Gamma} |_{V_\Gamma ^*}^2
	\label{y0}
\end{equation}
on the right hand side for all $t \in [0,T]$, 
with some positive constants $c_4$ and $c_5$. Notice that 
\eqref{yt} follows from \eqref{aux1} since we have 
\begin{align*}
	\bigl| \bar{u}_\Gamma (t) \bigr|_{V_\Gamma ^*}
	&
	= 
	\sup_{\substack{
	z_\Gamma \in V_\Gamma \\
	|z_\Gamma |_{V_\Gamma } \le 1}
	} \bigl| \langle 
	\bar{u}_\Gamma (t), z_\Gamma 
	\rangle _{V_\Gamma ^*, V_\Gamma }
	\bigr| 
	\nonumber \\
	& \le \sup_{\substack{
	z_\Gamma \in V_\Gamma \\
	|z_\Gamma |_{V_\Gamma } \le 1}
	} 
	\left|
	\int_{\Omega }^{} \nabla \bar{y}(t) \cdot \nabla {\mathcal R}z_\Gamma  dx + 
	\int_{\Gamma }^{} \nabla _\Gamma \bar{y}_\Gamma(t) \cdot \nabla _\Gamma z_\Gamma d\Gamma 
	\right|
	\nonumber \\
	& \le 
	\sup_{\substack{
	z_\Gamma \in V_\Gamma \\
	|z_\Gamma |_{V_\Gamma } \le 1}
	} 
	\left\{ c_2^{1/2} 
	{\revis 
	\bigl| \nabla \bar{y}(t) \bigr|_{H}
	}
	|z_\Gamma |_{V_\Gamma } + 
	{\revis 
	\bigl| \nabla _\Gamma \bar{y}_\Gamma (t) \bigr|_{H_\Gamma}
	}
	|z_\Gamma |_{V_\Gamma }
	\right\} 
\end{align*}
for all $t \in [0,T]$, where we have to use \eqref{co0} 
and the fact that $({\mathcal R}z_\Gamma,z_\Gamma) \in \boldsymbol{V}$ for all $z_\Gamma \in V_\Gamma$. 
Moreover, \eqref{y0} follows from \eqref{poin3} and \eqref{aux1} at time $0$, 
taking $\boldsymbol{z}=(\bar{y}(0), \bar{y}_\Gamma (0))$; 
of course, $c_5$ depends on $C_P$. 
Next, for all $t \in [0,T]$, we note that 
\begin{equation*}
	\int_{0}^{t} \bigl( \bar{\xi }(s), \bar{u}(s) \bigr)_H ds 
	\ge 0, 
	\quad 
	\int_{0}^{t} \bigl( \bar{\xi }_\Gamma (s), \bar{u}_\Gamma (s) \bigr)_{H_\Gamma } ds 
	\ge 0
\end{equation*}
by the monotonicity; 
\begin{equation*}
	\int_{0}^{t} \bigl( \pi \bigl( u^{(1)}(s) \bigr) -\pi \bigl( u^{(2)}(s) \bigr), \bar{u}(s) 
	\bigr)_H ds 
	\le L \int_{0}^{t} \bigl| \bar{u} (s) \bigr|_H^2 ds
\end{equation*}
and it will be treated by the {G}ronwall inequality;
\begin{align*}
	\int_{0}^{t} \bigl( \pi_\Gamma \bigl( u_\Gamma ^{(1)}(s) \bigr) 
	-\pi_\Gamma  \bigl( u_\Gamma ^{(2)}(s) \bigr), \bar{u}_\Gamma (s) 
	\bigr)_{H_\Gamma } ds
	& \le L_{\Gamma } 
	\int_{0}^{t} \bigl| \bar{u}_\Gamma  (s) \bigr|_{H_\Gamma }^2 ds 
	\\
	& \le \delta \int_{0}^{t}
	\bigl| \bar{u}_\Gamma  (s) \bigr|_{V_\Gamma }^2 ds 
	+ 
	C_\delta \int_{0}^{t}
	\bigl| \bar{u}_\Gamma  (s) \bigr|_{V_\Gamma^* }^2 ds 
\end{align*}
due to the compactness inequality \eqref{cp0}, 
also for this term we will use the the {G}ronwall inequality; 
the last two terms can be simply treated by the {Y}oung inequality; 
\begin{gather*}
	\int_{0}^{t} \bigl( \bar{f}(s), \bar{u}(s) \bigr)_H ds 
	\le \tilde{\delta} \int_{0}^{t} \bigl| \bar{u}(s) \bigr|_V^2 ds 
	+ \frac{1}{4\tilde{\delta} } \int_{0}^{t} \bigl| \bar{f}(s) \bigr|_{V^*}^2 ds,
	\\
	\int_{0}^{t} \bigl( \bar{f}_\Gamma (s), \bar{u}_\Gamma (s) \bigr)_{H_\Gamma } ds 
	\le \delta \int_{0}^{t} \bigl| \bar{u}_\Gamma (s) \bigr|_{V_\Gamma }^2 ds 
	+ \frac{1}{4\delta } \int_{0}^{t} \bigl| \bar{f}_\Gamma (s) \bigr|_{V_\Gamma ^*}^2 ds
\end{gather*}
for all $\delta, \tilde{\delta }>0$. 
Now we take advantage of \eqref{poin3}: indeed, $\bar{u}_\Gamma $ satisfies 
the zero mean value condition from \eqref{meanae}. 
Moreover, by \eqref{y0} we can 
recover from \eqref{last} the following inequality
\begin{align*}
	& \frac{\tau }{2} \bigl| \bar{u}(t) \bigr|_H^2 
	+ 
	c_4 \bigl| \bar{u}_\Gamma (t) \bigr|_{V_\Gamma ^*} ^2 
	+
	\frac{1}{2C_P} \int_{0}^{t} \bigl|\bar{u}(s) \bigr|_{V}^2 ds 
	+
	\frac{1}{2C_P} \int_{0}^{t} \bigl| \bar{u}_\Gamma (s) \bigr|_{V_\Gamma}^2 ds 
	\nonumber \\
	& \le 
	\frac{\tau }{2} | \bar{u}_0 |_H^2 
	+
	c_5 | \bar{u}_{0\Gamma } |_{V_\Gamma ^*}^2
	+ 
	 L \int_{0}^{t} \bigl| \bar{u} (s) \bigr|_H^2 ds
	+ 
	C_\delta \int_{0}^{t}
	\bigl| \bar{u}_\Gamma  (s) \bigr|_{V_\Gamma^* }^2 ds 
	\nonumber \\
	& \quad {}
	+ \frac{C_P}{2} \int_{0}^{t} \bigl| \bar{f}(s) \bigr|_{V^*}^2 ds
	+ C_P \int_{0}^{t} \bigl| \bar{f}_\Gamma (s) \bigr|_{V_\Gamma ^*}^2 ds 
\end{align*}
for all $t \in [0,T]$, that is, $\delta :=1/(4C_P)$ and $\tilde{\delta }:=1/(2C_P)$. 
Then, the continuous dependence \eqref{conti} follows from the application of 
the {G}ronwall inequality. \hfill $\Box$

\section*{Appendix}
\renewcommand{\theequation}{a.\arabic{equation}}
\setcounter{equation}{0}

We use the same notation as in the previous sections. 
We also use the following notation of function spaces. 
For each fixed $\varepsilon \in (0,1]$, 
define a subspace $\boldsymbol{H}_0^\varepsilon$ of $\boldsymbol{H}$ by 
$\boldsymbol{H}_0^\varepsilon:=\{ \boldsymbol{z} \in \boldsymbol{\boldsymbol{H}} : m_\varepsilon (\boldsymbol{z})=0 \}$, where 
$m_\varepsilon : \boldsymbol{H} \to \mathbb{R}$, 
\begin{equation*}
	m_\varepsilon (\boldsymbol{z}):= \frac{1}{\varepsilon |\Omega |+|\Gamma |}
	\left\{ \varepsilon \int_{\Omega }^{} z dx + \int_{\Gamma }^{} z_\Gamma d\Gamma \right\} 
	\quad \mbox{for all } \boldsymbol{z} \in \boldsymbol{H}.
\end{equation*}
Moreover, define an inner product of $\boldsymbol{H}$ by 
\begin{equation*}
	(\!( \boldsymbol{z}, \tilde{\boldsymbol{z}} )\!)_{\boldsymbol{H}}:=
	\varepsilon (z,\tilde{z})_H + (z_\Gamma, \tilde{z}_\Gamma  )_{H_\Gamma }
	\quad \mbox{for all } \boldsymbol{z} \in \boldsymbol{H}.
\end{equation*}
Then we see that the induced norm $\| \cdot  \|_{\boldsymbol{H}}$ and 
the standard norm $| \cdot |_{\boldsymbol{H}}$ are equivalent, because 
\begin{equation*}
	\| \boldsymbol{z} \|_{\boldsymbol{H}}^2 \le 
	| \boldsymbol{z} |_{\boldsymbol{H}}^2 \le 
	\frac{1}{\varepsilon } \| \boldsymbol{z} \|_{\boldsymbol{H}}^2
	\quad \mbox{for all } \boldsymbol{z} \in \boldsymbol{H}.
\end{equation*}
Next, we define $\boldsymbol{V}_0^\varepsilon 
:=\boldsymbol{V} \cap \boldsymbol{H}_0^\varepsilon $ with 
$| \boldsymbol{z}| _{\boldsymbol{V}_0^\varepsilon }:=\sqrt{a(\boldsymbol{z},\boldsymbol{z})}$ for all 
$\boldsymbol{z} \in \boldsymbol{V}_0^\varepsilon $, 
where we use a bilinear form 
$a(\cdot ,\cdot ):\boldsymbol{V} \times \boldsymbol{V} \to \mathbb{R}$ by 
\begin{equation*}
	a(\boldsymbol{z},\tilde{\boldsymbol{z}}):=
	\int_{\Omega }^{} \nabla z \cdot \nabla \tilde{z} dx 
	+\int_{\Gamma }^{} \nabla _\Gamma z_\Gamma \cdot \nabla _\Gamma \tilde{z}_\Gamma d\Gamma 
	\quad \mbox{for all } \boldsymbol{z},\tilde{\boldsymbol{z}} \in \boldsymbol{V}.
\end{equation*} 
Let us define a linear bounded operator 
$\boldsymbol{F}: \boldsymbol{V}_0^\varepsilon \to (\boldsymbol{V}_0^\varepsilon )^*$ by 
\begin{equation*}
	\langle \boldsymbol{F}\boldsymbol{z}, \tilde{\boldsymbol{z}} \rangle_{(\boldsymbol{V}_0^\varepsilon )^*, \boldsymbol{V}_0^\varepsilon }
	:=a(\boldsymbol{z},\tilde{\boldsymbol{z}})
	\quad \mbox{for all } \boldsymbol{z},\tilde{\boldsymbol{z}} \in \boldsymbol{V}_0^\varepsilon.
\end{equation*}
Then, there exists a positive constant $c_P^{(\varepsilon )}$ such that 
\begin{equation*}
	\| \boldsymbol{z} \|_{\boldsymbol{H}}^2 
	\le 
	| \boldsymbol{z} |_{\boldsymbol{H}}^2 \le 
	c_P^{(\varepsilon )} | \boldsymbol{z} |_{\boldsymbol{V}_0^\varepsilon }^2
	\quad \mbox{for all } \boldsymbol{z} \in \boldsymbol{V}_0^\varepsilon,
\end{equation*}
see e.g., \cite[Appendix]{CF15b}. Thus, we see that 
$| \cdot |_{\boldsymbol{V}_0^\varepsilon }$ and 
standard $| \cdot |_{\boldsymbol{V}}$ are equivalent norm of 
$\boldsymbol{V}_{0}^\varepsilon $ and 
$\boldsymbol{F}$ is the duality mapping from 
$\boldsymbol{V}_0^\varepsilon $ to $(\boldsymbol{V}_0^\varepsilon )^*$.  
We also define the inner product in $(\boldsymbol{V}_0^\varepsilon )^*$ by
\begin{equation*}
	(\boldsymbol{z}_1^*,\boldsymbol{z}_2^*)_{(\boldsymbol{V}_0^\varepsilon )^*}
	:=
	\langle \boldsymbol{z}_1^*, 
	\boldsymbol{F}^{-1}\boldsymbol{z}_2^* \rangle_{(\boldsymbol{V}_0^\varepsilon )^*, \boldsymbol{V}_0^\varepsilon }
	\quad \mbox{for all } \boldsymbol{z}_1^*,
	\boldsymbol{z}_2^* \in (\boldsymbol{V}_0^\varepsilon)^*.
\end{equation*}
Thanks to \cite[Appendix]{CF15b} again, we obtain $\boldsymbol{V}_0^\varepsilon 
\mathop{\hookrightarrow} \mathop{\hookrightarrow}
\boldsymbol{H}_0^\varepsilon 
\mathop{\hookrightarrow} \mathop{\hookrightarrow}
(\boldsymbol{V}_0^\varepsilon )^*$.

\paragraph{Lemma A.} 
{\it An operator $\boldsymbol{P}_\varepsilon :\boldsymbol{H} \to \boldsymbol{H}_0^\varepsilon $ defined by 
\begin{equation*}
	\boldsymbol{P}_\varepsilon \boldsymbol{z} :=\boldsymbol{z}-m_\varepsilon (\boldsymbol{z}) \boldsymbol{1}
	\quad \mbox{for all } \boldsymbol{z} \in \boldsymbol{H}
\end{equation*}
is the projection from $\boldsymbol{H}$ to $\boldsymbol{H}_0^\varepsilon$ with 
respect to $\| \cdot \|_{\boldsymbol{H}}$-norm, namely}
\begin{equation}
	\| \boldsymbol{z}-\boldsymbol{P}_\varepsilon \boldsymbol{z} \|_{\boldsymbol{H}} 
	\le
	\| \boldsymbol{z} - \boldsymbol{y} \|_{\boldsymbol{H}}
	\quad \mbox{\it for all } \boldsymbol{y} \in \boldsymbol{H}_0^\varepsilon.
	\label{defpr}
\end{equation}

\paragraph{Proof.} Firstly, we see that 
\begin{align*}
	m_\varepsilon (\boldsymbol{P}_\varepsilon \boldsymbol{z}) & = 
	 \frac{1}{\varepsilon |\Omega |+|\Gamma |}
	\left\{ \varepsilon \int_{\Omega }^{} \bigl( 
	z-m_\varepsilon (\boldsymbol{z}) 
	\bigr) dx + \int_{\Gamma }^{} 
	\bigl( 
	z_\Gamma -m_\varepsilon (\boldsymbol{z}) 
	\bigr) 
	d\Gamma \right\} 
	\nonumber \\
	& = m_\varepsilon ( \boldsymbol{z}) 
	- \frac{m_\varepsilon ( \boldsymbol{z}) }{\varepsilon |\Omega |+|\Gamma |}
	\left\{ \varepsilon |\Omega |+|\Gamma | \right\}
	\nonumber \\
	& = 0
	\quad \mbox{for all }\boldsymbol{z} \in \boldsymbol{H}.
\end{align*}
Next, we see that 
\begin{align*}
	& (\!( 
	\boldsymbol{z}-\boldsymbol{P}_\varepsilon \boldsymbol{z}, 
	\boldsymbol{y}-\boldsymbol{P}_\varepsilon \boldsymbol{z}
	)\!)_{\boldsymbol{H}} 
	\\
	& \quad =  \bigl( \! \bigl( 
	m_\varepsilon (\boldsymbol{z}) \boldsymbol{1},  
	\boldsymbol{y}-\boldsymbol{z}+ m_\varepsilon  (\boldsymbol{z}) \boldsymbol{1}
	\bigr)\! \bigr)_{\boldsymbol{H}} 
	\\
	& \quad = m_\varepsilon (\boldsymbol{z}) 
	\left\{ \varepsilon \int_{\Omega }^{} \bigl( 
	y-z
	\bigr) dx + \int_{\Gamma }^{} 
	\bigl( 
	y_\Gamma -z_\Gamma 
	\bigr) 
	d\Gamma \right\} 
	+ m_\varepsilon (\boldsymbol{z})^2 
	\bigl\{ \varepsilon |\Omega | + |\Gamma | \bigr\} 
	\\
	& \quad = m_\varepsilon (\boldsymbol{z})  m_\varepsilon (\boldsymbol{y})
	\bigl\{ \varepsilon |\Omega | + |\Gamma | \bigr\} 
	\\ 
	& \quad =0 \quad \mbox{for all } \boldsymbol{y} \in \boldsymbol{H}_0^\varepsilon. 
\end{align*}
Therefore
\begin{align*}
	\| \boldsymbol{z}-\boldsymbol{P}_\varepsilon \boldsymbol{z} \|_{\boldsymbol{H}}^2 
	& = (\!( 
	\boldsymbol{z}-\boldsymbol{P}_\varepsilon \boldsymbol{z}, 
	\boldsymbol{z}-\boldsymbol{y}+\boldsymbol{y}-\boldsymbol{P}_\varepsilon \boldsymbol{z}
	)\!)_{\boldsymbol{H}} \\
	& = (\!( 
	\boldsymbol{z}-\boldsymbol{P}_\varepsilon \boldsymbol{z}, 
	\boldsymbol{z}-\boldsymbol{y}
	)\!)_{\boldsymbol{H}} \\
	& \le \frac{1}{2} \| 
	\boldsymbol{z}-\boldsymbol{P}_\varepsilon \boldsymbol{z}
	\|_{\boldsymbol{H}}^2 
	+ \frac{1}{2}
	\| \boldsymbol{z}-\boldsymbol{y} \|_{\boldsymbol{H}}^2,
\end{align*}
namely \eqref{defpr} holds. \hfill $\Box$ 

\bigskip
We also easily obtain the following conditions:
\begin{gather*}
	( \! ( \boldsymbol{z}^*, \boldsymbol{P}_\varepsilon \boldsymbol{z} 
	) \! )_{\boldsymbol{H}}
	= (\!( \boldsymbol{z}^*, \boldsymbol{z} )\!)_{\boldsymbol{H}}
	\quad \mbox{for all }\boldsymbol{z}^* \in \boldsymbol{H}_0^\varepsilon 
	\ \mbox{and } \boldsymbol{z} \in \boldsymbol{H}, \\
	|\boldsymbol{P}_\varepsilon \boldsymbol{z}|_{\boldsymbol{V}_0^\varepsilon }
	\le |\boldsymbol{z}|_{\boldsymbol{V}}
	\quad \mbox{for all } \boldsymbol{z} \in \boldsymbol{V}. 
\end{gather*}
Incidentally, we note that another possibility of projection is given by
\begin{equation*}
	\tilde{\boldsymbol{P}}_\varepsilon \boldsymbol{z}
	:= \left( z - \frac{1}{|\Omega |}
	\int_{\Omega }^{} z dx,
	z_\Gamma  - \frac{1}{|\Gamma |}\int_{\Gamma }^{} z_\Gamma d\Gamma \right)
	\quad \mbox{for all } \boldsymbol{z} \in \boldsymbol{H},
\end{equation*}
and $\tilde{\boldsymbol{P}}_\varepsilon $ 
is actually the projection from $\boldsymbol{H}$ to $\boldsymbol{H}_0^\varepsilon$ with 
respect to the standard norm 
(cf.\ \eqref{defpr}). However, this choice is not suitable from the viewpoint of the trace condition. 
Indeed, $\tilde{\boldsymbol{P}}_\varepsilon \boldsymbol{z} \notin \boldsymbol{V}$ 
even if $\boldsymbol{z} \in \boldsymbol{V}$.

\medskip
Next, we prepare suitable 
approximations for $\boldsymbol{f}$ and $\boldsymbol{u}_0$, 
which satisfy assumptions (A1) and (A2), respectively.

\paragraph{Lemma B.}
{\it For each $\varepsilon \in (0,1]$, the solution of the following {C}auchy problem 
\begin{numcases}
	{} 
	\varepsilon \boldsymbol{f}_\varepsilon '(s) + \boldsymbol{f}_\varepsilon(s) =\boldsymbol{f}(s) 
	\quad \mbox{in } \boldsymbol{H}, \quad \mbox{for a.a.\ } s \in (0,T), 
	\label{apfex}\\
	\boldsymbol{f}_\varepsilon (0) = \boldsymbol{0}
	\quad \mbox{in } \boldsymbol{H} \nonumber 
\end{numcases}
satisfies \eqref{apf1} and \eqref{apf2}. }

\paragraph{Proof.}
Note that $\boldsymbol{f}_\varepsilon (0)=\boldsymbol{0} \in \boldsymbol{V}$. 
Multiplying the above equation by the solution 
$\boldsymbol{f}_\varepsilon :=(f_\varepsilon , f_{\Gamma, \varepsilon }) 
\in H^1(0,T;\boldsymbol{H})$, 
integrating it over $[0,t]$ 
and using {Y}oung inequality we have 
\begin{equation}
	\frac{\varepsilon }{2} \bigl| \boldsymbol{f}_\varepsilon (t) \bigr|_{\boldsymbol{H}}^2
	+ \frac{1}{2}\int_{0}^{t} \bigl| \boldsymbol{f}_\varepsilon (s) \bigr|_{\boldsymbol{H}}^2 ds
	\le \frac{1}{2}\int_{0}^{t} \bigl| \boldsymbol{f} (s) \bigr|_{\boldsymbol{H}}^2 ds
	\label{apfbound}
\end{equation}
for all $t \in [0,T]$. 
Thus, 
there exist a subsequence of $\varepsilon $ (not relabeled) and some 
limit function $\tilde{\boldsymbol{f}} \in L^2(0,T;\boldsymbol{H})$ such that 
\begin{equation*}
	\boldsymbol{f}_\varepsilon \to 
	\tilde{\boldsymbol{f}} \quad \mbox{weakly in } L^2(0,T;\boldsymbol{H})
\end{equation*}
{\takeshi as $\varepsilon \searrow 0$}. Moreover, from \eqref{apfex} and \eqref{apfbound}, there exists a positive constant $M_{10}$ such that 
\begin{gather*}
	|\varepsilon \boldsymbol{f}_\varepsilon ' |_{L^2(0,T;\boldsymbol{H})} = 
	|\boldsymbol{f}-\boldsymbol{f}_\varepsilon |_{L^2(0,T;\boldsymbol{H})} \le M_{10},
	\\
	\varepsilon ^{1/2} |\boldsymbol{f}_\varepsilon|_{L^\infty (0,T;\boldsymbol{H})}\le M_{10},
\end{gather*}
that is, there 
exist a subsequence of $\varepsilon $ (not relabeled) and some 
limit function $\dot{\boldsymbol{f}} \in L^2(0,T;\boldsymbol{H})$ such that 
\begin{gather*}
	\varepsilon \boldsymbol{f}'_\varepsilon \to 
	\dot{\boldsymbol{f}} \quad \mbox{weakly in } L^2(0,T;\boldsymbol{H}), 
	\nonumber \\
	\varepsilon \boldsymbol{f}_\varepsilon \to 
	\boldsymbol{0} \quad \mbox{strongly in } {\revis C \bigl( [0,T];\boldsymbol{H} \bigr)}
\end{gather*}
{\takeshi as $\varepsilon \searrow 0$}, that is $\dot{\boldsymbol{f}}=\boldsymbol{0}$. 
From \eqref{apfex}, this implies that 
$\tilde{\boldsymbol{f}}=\boldsymbol{f}$. Therefore, by using \eqref{apfbound} again 
\begin{equation*}
	\limsup_{\varepsilon \to 0} |\boldsymbol{f}_\varepsilon |_{L^2(0,T;\boldsymbol{H})} 
	\le |\boldsymbol{f}|_{L^2(0,T;\boldsymbol{H})},
\end{equation*}
which gives us the {\colli convergence of norms and consequently the} strong convergence \eqref{apf1}. 
Next we show the required order of the convergence \eqref{apf2}. Recall the equation for 
$f_{\Gamma, \varepsilon} \in H^1(0,T;H_\Gamma)$ as follows:
\begin{numcases}
	{} 
	\varepsilon f_{\Gamma, \varepsilon }'(s) + f_{\Gamma, \varepsilon}(s) = f_\Gamma (s) 
	\quad \mbox{in } H_\Gamma, \quad \mbox{for a.a.\ } s \in (0,T), 
	\label{apfgex}\\
	f_{\Gamma, \varepsilon} (0) = 0
	\quad \mbox{in } H_\Gamma , \nonumber 
\end{numcases}
From \eqref{apfgex} we see that 
\begin{equation*}
	|f_\Gamma -f_{\Gamma ,\varepsilon }|_{L^2(0,T;H_\Gamma )} 
	= \varepsilon |f_{\Gamma , \varepsilon }'|_{L^2(0,T;H_\Gamma )},
\end{equation*}
therefore it is enough to prove that there exists a positive constant $C_0$ such that
\begin{equation*}
	\varepsilon ^{1/2}|f_{\Gamma , \varepsilon }'|_{L^2(0,T;H_\Gamma )} \le C_0.
\end{equation*}
Test \eqref{apfgex} by $f'_{\Gamma, \varepsilon }(s)$, 
integrating the resultant with respect to 
time, we obtain 
\begin{align*}
	\varepsilon \int_{0}^{t} \bigl| f_{\Gamma, \varepsilon} ' (s) \bigr|_{H_\Gamma }^2 ds 
	+ \frac{1}{2} \bigl| f_{\Gamma, \varepsilon  }(t) \bigr|^2_{H_\Gamma } 
	& = \bigl( f_{\Gamma }(t), f_{\Gamma, \varepsilon  }(t) \bigr)_{H_\Gamma }
	- \int_{0}^{t} \bigl( f_{\Gamma }'(s), f_{\Gamma, \varepsilon  }(s) \bigr)_{H_\Gamma } ds 
	\\
	& \le \bigl| f_{\Gamma }(t) \bigr|_{H_\Gamma }^2 
	 + \frac{1}{4} \bigl| f_{\Gamma, \varepsilon  }(t) \bigr|_{H_\Gamma }^2
	+ \int_{0}^{t} \bigl| f_{\Gamma }'(s) \bigr|_{H_\Gamma } 
	\bigl| f_{\Gamma, \varepsilon  }(s) \bigr|_{H_\Gamma } ds 
\end{align*}
for all $t \in [0,T]$. Therefore, using {G}ronwall inequality we see that 
there exists a positive constant $C_0$, depending on 
$|f_\Gamma |_{C([0,T];H_\Gamma )}$ and $|f_\Gamma '|_{L^1(0,T;H_\Gamma )}$, 
independent of $\varepsilon \in (0,1]$, 
such that
\begin{equation*}
	\varepsilon ^{1/2} | f_{\Gamma, \varepsilon} ' |_{L^2(0,T;H_\Gamma) } 
	+ | f_{\Gamma, \varepsilon  }|_{L^\infty (0,T;H_\Gamma )}
	\le C_0 
\end{equation*}
for all $\varepsilon \in (0,1]$. \hfill $\Box$

\paragraph{Lemma C.}
{\it For each $\varepsilon \in (0,1]$, the solution 
$\boldsymbol{u}_{0,\varepsilon }
:=(u_{0,\varepsilon }, u_{0\Gamma ,\varepsilon }) \in \boldsymbol{W} \cap \boldsymbol{V} $
of the following elliptic system
\begin{numcases}
	{} 
	u_{0,\varepsilon }-\varepsilon \Delta u_{0, \varepsilon } =u_0
	\quad \mbox{a.e.\ in } \Omega, 
	\label{apuex}\\
	(u_{0, \varepsilon})_{|_\Gamma } 
	= u_{0\Gamma, \varepsilon  }, \quad u_{0\Gamma,\varepsilon }
	+ \varepsilon \partial _{\boldsymbol{\nu }} u_{0, \varepsilon } 
	- \varepsilon \Delta_\Gamma  u_{0\Gamma, \varepsilon } =u_{0\Gamma }
	\quad \mbox{a.e.\ on } \Gamma 
	\label{apugex}
\end{numcases}
satisfies \eqref{apu} and 
\eqref{apbound} with the required regularity
$( -\Delta u_{0, \varepsilon }, 
\partial _{\boldsymbol{\nu }} u_{0, \varepsilon } - \Delta _\Gamma u_{0\Gamma, \varepsilon} ) 
\in \boldsymbol{V}$. }

\paragraph{Proof.} The strategy of the proof is similar to 
one of \cite[Lemma~A.1]{CF16}. 
By virtue of \cite[Lemma~C]{CF15b}, 
there exists $\boldsymbol{u}_{0,\varepsilon }
:=(u_{0,\varepsilon }, u_{0\Gamma ,\varepsilon }) \in \boldsymbol{W} \cap \boldsymbol{V} $ such that 
$\boldsymbol{u}_{0,\varepsilon }$ satisfies the 
system \eqref{apuex}--\eqref{apugex}. Moreover 
the assumption (A2) gives us the required regularity
$( -\Delta u_{0, \varepsilon }, 
\partial _{\boldsymbol{\nu }} u_{0, \varepsilon } - \Delta _\Gamma u_{0\Gamma, \varepsilon} ) 
\in \boldsymbol{V}$.
Next we show \eqref{apu}. 
Indeed, testing the first equation \eqref{apuex} by $u_{0, \varepsilon}$, 
the second equation \eqref{apugex} by $u_{0\Gamma, \varepsilon}$, 
adding them, and using the {Y}oung inequality 
we obtain 
\begin{gather*}
	\frac{1}{2}\int_{\Omega }^{} |u_{0, \varepsilon } |^2dx 
	+ \frac{1}{2}\int_{\Gamma }^{} |u_{0\Gamma, \varepsilon  }|^2d\Gamma 
	+ \varepsilon \int_{\Omega }^{} |\nabla u_{0, \varepsilon }|^2dx 
	+ \varepsilon \int_{\Gamma }^{} |\nabla_\Gamma u_{0\Gamma, \varepsilon }|^2d\Gamma  
	\\
	\le \frac{1}{2}\int_{\Omega }^{} |u_{0} |^2dx 
	+ \frac{1}{2}\int_{\Gamma }^{} |u_{0\Gamma}|^2d\Gamma. 
\end{gather*}
Therefore, 
$\{ u_{0, \varepsilon} \}_{\varepsilon \in (0,1]}$ is bounded in $H$,  
$\{ u_{0\Gamma, \varepsilon} \}_{\varepsilon \in (0,1]}$ is bounded in $H_\Gamma $,  
$\{ \varepsilon^{1/2} \nabla u_{0, \varepsilon} \}_{\varepsilon \in (0,1]}$ is bounded in {\revis $H$}, 
and 
$\{ \varepsilon^{1/2} \nabla_\Gamma  u_{0, \varepsilon} \}_{\varepsilon \in (0,1]}$ 
is bounded in {\revis $H_{\Gamma}$}, respectively. 
These give us (see \eqref{apuex}--\eqref{apugex}) 
\begin{gather*}
	u_{0, \varepsilon} \to u_0 \quad \mbox{weakly in } H, \quad 
	\varepsilon u_{0, \varepsilon} \to 0 \quad \mbox{strongly in } V, \\
	u_{0\Gamma , \varepsilon} \to u_{0\Gamma } \quad \mbox{weakly in } H_\Gamma, \quad 
	\varepsilon u_{0\Gamma , \varepsilon} \to 0 \quad \mbox{strongly in } V_\Gamma
\end{gather*}
{\takeshi as $\varepsilon \searrow 0$}. Moreover, we have
\begin{equation*}
	\limsup_{\varepsilon \to 0 } \left( \int_{\Omega }^{} |u_{0, \varepsilon } |^2dx 
	+\int_{\Gamma }^{} |u_{0\Gamma, \varepsilon  }|^2d\Gamma  \right) 
	\le \int_{\Omega }^{} |u_{0} |^2dx 
	+ \int_{\Gamma }^{} |u_{0\Gamma  } |^2d\Gamma, 
\end{equation*}
which {\revis entails} that 
\begin{equation*}
	\boldsymbol{u}_{0,\varepsilon }
	\to 
	\boldsymbol{u}_0
	\quad \mbox{strongly in } \boldsymbol{H}
\end{equation*} 
{\takeshi as $\varepsilon \searrow 0$}. 
Next, from the definition of the subdifferential 
with \eqref{apuex} and \eqref{apugex} we see that 
\begin{align*}
	\int_{\Omega }^{}
	\widehat{\beta }_\lambda (u_{0, \varepsilon }) dx 
	-
	\int_{\Omega }^{} 
	\widehat{\beta }_\lambda (u_{0}) dx 
	& \le 
	\int_{\Omega }^{} 
	(u_{0, \varepsilon }-u_0) \beta _\lambda (u_{0, \varepsilon }) dx 
	\\
	& = -\int_{\Omega }^{} \varepsilon \beta '_\lambda (u_{0,\varepsilon } )
	|\nabla u_{0, \varepsilon }|^2dx 
	+ \varepsilon 
	\int_{\Gamma }^{} 
	\partial _{\boldsymbol{\nu }} u_{0, \varepsilon }
	\beta _\lambda (u_{0\Gamma , \varepsilon }) d\Gamma \\
	& \le 
	\int_{\Gamma }^{} (u_{0\Gamma }-u_{0\Gamma, \varepsilon })
	\beta _\lambda (u_{0\Gamma, \varepsilon}) d\Gamma 
	- \varepsilon \int_{\Gamma }^{} \beta _\lambda '(u_{0\Gamma, \varepsilon })
	|\nabla u_{0\Gamma, \varepsilon  }|^2 d\Gamma 
	\\
	& \le \int_{\Gamma }^{} \widehat{\beta }_\lambda (u_{0\Gamma }) d\Gamma 
	- \int_{\Gamma }^{} \widehat{\beta }_\lambda (u_{0\Gamma, \varepsilon }) d\Gamma, 
\end{align*}
that is, 
\begin{equation}
	\int_{\Omega }^{}
	\widehat{\beta }_\lambda (u_{0, \varepsilon }) dx 
	+
	\int_{\Gamma }^{} \widehat{\beta }_\lambda (u_{0\Gamma, \varepsilon }) d\Gamma
	\le 
	\int_{\Omega }^{} 
	\widehat{\beta }(u_{0}) dx 
	+
	\int_{\Gamma }^{} \widehat{\beta }(u_{0\Gamma }) d\Gamma, 
\end{equation}
for all $\varepsilon \in (0,1]$, where \eqref{prim} has been used. 
Now, due to \eqref{A5e} we obtain that 
\begin{equation*}
	\widehat{\beta }_\lambda (r) = \int_{0}^{r} \beta _\lambda (s) ds \le 
	\rho \int_{0}^{r} \beta _{\Gamma, \lambda} (s) ds +c_0 |r| 
	= \rho \widehat{\beta }_{\Gamma, \lambda} (r) + c_0 |r|
\end{equation*} 
for all $r \in \mathbb{R}$. Indeed, in the case $r \ge 0$, because of the fact that 
$\beta _\lambda (0)=\beta _{\Gamma,\lambda }(0)=0$, we see from \eqref{A5e} that 
\begin{equation*}
	\int_{0}^{r} \beta _\lambda (s) ds = 
	\int_{0}^{r} \bigl| \beta _\lambda (s)\bigr| ds\le 
	\rho \int_{0}^{r} \bigl| \beta _{\Gamma, \lambda} (s) \bigr| ds +c_0 r 
	= \rho \int_{0}^{r} \beta _{\Gamma, \lambda} (s) ds +c_0 |r|. 
\end{equation*}
In the case $r < 0$, we have
\begin{equation*}
	\int_{0}^{r} \beta _\lambda (s) ds = 
	\int_{r}^{0} \bigl| \beta _\lambda (s)\bigr| ds\le 
	\rho \int_{r}^{0} \bigl| \beta _{\Gamma, \lambda} (s) \bigr| ds  - c_0 r 
	= \rho \int_{0}^{r} \beta _{\Gamma, \lambda} (s) ds +c_0 |r|. 
\end{equation*}
Then, it turns out that 
\begin{align*}
	\int_{\Gamma }^{} \widehat{\beta }_\lambda (u_{0\Gamma }) d\Gamma 
	& \le 
	\rho \int_{\Gamma }^{} \widehat{\beta }_{\Gamma, \lambda} (u_{0\Gamma }) d\Gamma 
	 + c_0 |u_{0\Gamma }|_{L^1(\Gamma )} \\
	& \le \rho \int_{\Gamma }^{} \widehat{\beta }_{\Gamma} (u_{0\Gamma }) d\Gamma 
	 + c_0 |u_{0\Gamma }|_{L^1(\Gamma )}. 
\end{align*}
Therefore, $\widehat{\beta}(u_{0\Gamma }) \in L^1(\Gamma )$ due to 
the {F}atou lemma and the almost everywhere convergence of $\widehat{\beta }_{\lambda }(u_{0\Gamma })$ to 
$\widehat{\beta }(u_{0\Gamma })$. 
Concerning our approximation, testing \eqref{apuex} by $-\Delta u_{0,\varepsilon }$ and using 
\eqref{apugex}, we get 
\begin{gather*}
	\int_{\Omega }^{} \nabla (u_{0,\varepsilon }-u_0) \cdot \nabla u_{0\varepsilon } dx 
	+ \frac{1}{\varepsilon } \int_{\Gamma }^{}
	|u_{0\Gamma ,\varepsilon }-u_{0\Gamma }|^2 d\Gamma 
	+ \varepsilon \int_{\Omega }^{} |\Delta u_{0\varepsilon }|^2 dx 
	\\
	{}
	+ \int_{\Gamma }^{} \nabla _\Gamma u_{0\Gamma, \varepsilon} \cdot 
	\nabla _\Gamma (u_{0\Gamma ,\varepsilon }-u_{0\Gamma }) d\Gamma =0.
\end{gather*}
Therefore, using the {Y}oung inequality we deduce 
\begin{gather*}
	\frac{1}{2} \int_{\Omega }^{} | \nabla u_{0,\varepsilon }|^2dx 
	+ \frac{1}{2} \int_{\Gamma }^{} |\nabla _\Gamma u_{0\Gamma, \varepsilon} |^2 d\Gamma 
	+ \frac{1}{\varepsilon } \int_{\Gamma }^{}
	|u_{0\Gamma ,\varepsilon }-u_{0\Gamma }|^2 d\Gamma 
	+ \varepsilon \int_{\Omega }^{} |\Delta u_{0\varepsilon }|^2 dx 
	\\
	\le \frac{1}{2}|\nabla u_0|^2 dx 
	+
	\frac{1}{2}
	\int_{\Gamma }^{} |\nabla _\Gamma u_{0\Gamma}|^2 d\Gamma. 
\end{gather*}
Thus, we obtain that there exists a positive constant $\tilde{C}_0$ such that
\begin{gather*}
	|u_{0\Gamma ,\varepsilon }-u_{0\Gamma }|_{H_\Gamma } \le  \varepsilon ^{1/2}\tilde{C}_0,
	\\
	\varepsilon ^{1/2} |\Delta u_{0\varepsilon }|_{H} \le \tilde{C}_0 \quad \mbox{for all } 
	\varepsilon \in (0,1], 
	\\
	\boldsymbol{u}_{0, \varepsilon } \to 
	\boldsymbol{u}_0 
	\quad \mbox{strongly in } \boldsymbol{V}
\end{gather*}
{\takeshi as $\varepsilon \searrow 0$, because, 
$\boldsymbol{u}_{0, \varepsilon } \to \boldsymbol{u}_0$ 
weakly in $\boldsymbol{V}$ and 
$|\boldsymbol{u}_{0, \varepsilon }|_{\boldsymbol{V}} \to |\boldsymbol{u}_0|_{\boldsymbol{V}}$ as $\varepsilon \searrow 0$.}
Then, from \eqref{apugex} we can also infer that 
\begin{equation*} 
	\varepsilon ^{1/2}
	|\partial _{\boldsymbol{\nu }} 
	u_{0,\varepsilon }-\Delta_\Gamma  u_{0\Gamma, \varepsilon}|_{H_\Gamma } \le \tilde{C}_0.
\end{equation*}
Now, if we test \eqref{apugex} by $\beta _{\Gamma,\lambda}(u_{0\Gamma, \varepsilon})$, then we obtain 
\begin{align*}
	\int_{\Gamma }^{} \widehat{\beta }_{\Gamma,\lambda} (u_{0\Gamma, \varepsilon}) d\Gamma 
	-\int_{\Gamma }^{} \widehat{\beta }_{\Gamma,\lambda} (u_{0\Gamma}) d\Gamma 
	& \le \int_{\Gamma }^{} (-\varepsilon \partial _{\boldsymbol{\nu }} 
	u_{0,\varepsilon }+ \varepsilon \Delta_\Gamma u_{0\Gamma, \varepsilon}) 
	\beta _{\Gamma,\lambda} (u_{0\Gamma, \varepsilon})
	d\Gamma \\
	& \le \tilde{C}_0 \varepsilon ^{1/2} \bigl| 
	\beta _{\Gamma,\lambda} (u_{0\Gamma, \varepsilon}) 
	\bigr|_{H_\Gamma }\\
	& \le \frac{\tilde{C}_0 \varepsilon ^{1/2}}{\lambda }|u_{0\Gamma, \varepsilon}|_{H_\Gamma },
\end{align*}
hence from \eqref{prim} there exists a positive constant $\bar{C}_0$ independent of $\varepsilon, \lambda \in (0,1]$, such that
\begin{equation*}
	\int_{\Gamma }^{} \widehat{\beta }_{\Gamma,\lambda} (u_{0\Gamma, \varepsilon}) d\Gamma 
	\le \int_{\Gamma }^{} \widehat{\beta }_{\Gamma,\lambda} (u_{0\Gamma}) d\Gamma 
	+ \frac{\bar{C}_0 \varepsilon ^{1/2}}{\lambda }
	\le \int_{\Gamma }^{} \widehat{\beta }_{\Gamma} (u_{0\Gamma}) d\Gamma 
	+ \frac{\bar{C}_0 \varepsilon ^{1/2}}{\lambda }.
\end{equation*} 
Thus, we get the conclusion. \hfill $\Box$

\paragraph{Proof of Proposition 3.1.}
For each $\varepsilon , \lambda \in (0,1]$, 
let us consider the 
{C}auchy problem for the following equivalent evolution equation: 
\begin{gather}
	\boldsymbol{v}'(t)+ \boldsymbol{F} 
	\bigl( \boldsymbol{P}_\varepsilon \boldsymbol{\mu }(t) \bigr)
	=\boldsymbol{0}
	\quad \mbox{in } ( \boldsymbol{V}_0^\varepsilon )^*, 
	\ \mbox{for a.a.\ }  t\in (0,T), 
	\label{Ap2}\\
	\boldsymbol{M} \boldsymbol{\mu }(t)
	= \boldsymbol{T} \boldsymbol{v}'(t)
	+ \boldsymbol{M} \partial \varphi \bigl( \boldsymbol{v}(t) \bigr)
	+ \boldsymbol{\beta }_\lambda \bigl( \boldsymbol{u}(t) \bigr)
	+ \boldsymbol{\pi }\bigl( \boldsymbol{u}(t) \bigr)
	-\boldsymbol{f}_\varepsilon (t)
	\nonumber \\
	\quad \mbox{in } \boldsymbol{H}, 
	\ \mbox{for a.a.\ }  t\in (0,T), 
	\label{Ap3}\\
	\boldsymbol{u}(t) = \boldsymbol{v}(t)+m_\varepsilon (\boldsymbol{u}_{0, \varepsilon })\boldsymbol{1}, 
	\quad 
	\boldsymbol{v}(0) = \boldsymbol{v}_0:=\boldsymbol{u}_{0, \varepsilon}-m_\varepsilon (\boldsymbol{u}_{0, \varepsilon})
	\quad \mbox{in }\boldsymbol{H}_0^\varepsilon, 
	\label{Ap4}
\end{gather}
(cf.\ \cite[Section~2.3]{CF15b}), where we define 
\begin{equation*}
	\boldsymbol{M}:=
	\begin{pmatrix}
	\varepsilon & 0 \\
	0 & 1
	\end{pmatrix}, 
	\quad 
	\boldsymbol{T}:=
	\begin{pmatrix}
	\tau & 0 \\
	0 & \varepsilon 
	\end{pmatrix}, 
\end{equation*}
and $\varphi :\boldsymbol{H}_0 ^\varepsilon \to [0,+\infty ]$ by 
\begin{equation*}
	\varphi (\boldsymbol{z}) 
	:=
	\begin{cases}
	\displaystyle 
	\frac{1}{2}
	\int_{\Omega }^{}|\nabla z|^2dx 
	+
	\frac{1}{2}
	\int_{\Gamma }^{} |\nabla _\Gamma z_\Gamma |^2d\Gamma & \mbox{if } \boldsymbol{z} \in \boldsymbol{V}_0^\varepsilon, \\
	+ \infty & \mbox{if } \boldsymbol{z} \in \boldsymbol{H}_0^\varepsilon \setminus \boldsymbol{V}_0^\varepsilon. 
	\end{cases}
\end{equation*}
Then, we see that $\varphi $ is proper, lower semicontinuous and convex on $\boldsymbol{H}_0^\varepsilon $ 
and the subdifferential $\partial  \varphi $ on $\boldsymbol{H}_0^\varepsilon $ is characterized by 
$\partial \varphi (\boldsymbol{z})=(-(1/\varepsilon )\Delta z, \partial _{\boldsymbol{\nu }}z-\Delta _\Gamma z_\Gamma )$ with 
$\boldsymbol{z} \in D(\partial \varphi )=\boldsymbol{W}\cap \boldsymbol{V}_0^\varepsilon $. 
The {C}auchy problem \eqref{Ap2}--\eqref{Ap4} can be solved by applying 
\cite[Sections~4.2--4.4]{CF15b}, based on the abstract theory of doubly nonlinear evolution equation \cite{CV90}, because 
all assumptions for $\boldsymbol{f}_\varepsilon $ 
and $\boldsymbol{u}_{0,\varepsilon }$
in order to  
obtain the strong solution are satisfied. \hfill $\Box$

\section*{Acknowledgments}

The present paper benefits from the support of the MIUR-PRIN Grant 2015PA5MP7 
``Calculus of Variations'', the GNAMPA 
(Gruppo Nazionale per l'Analisi Matematica, la Probabilit\`a e le loro Applicazioni) 
of INdAM (Istituto Nazionale di Alta Matematica) 
and the IMATI -- C.N.R. Pavia for PC; 
the support of the JSPS KAKENHI 
Grant-in-Aid for Scientific Research(C), Grant Number 17K05321 for TF.

\end{document}